\documentclass[twoside, 12pt]{article}
\usepackage[latin1]{inputenc}
\usepackage{amscd}
\usepackage{amsmath}
\usepackage{amsfonts}
\usepackage{amssymb}
\usepackage{amsthm}
\usepackage{comment}
\usepackage[francais,english]{babel}
\usepackage[matrix,arrow]{xy}

\author{S. Gouëzel}

\newtheorem{thm}{Theorem}[section]
\newtheorem*{thm*}{Theorem}
\newtheorem{prop}[thm]{Proposition}
\newtheorem{lem}[thm]{Lemma}
\newtheorem{cor}[thm]{Corollary}
\newtheorem{defn}[thm]{Definition}

\theoremstyle{definition}

\newtheorem*{rmq}{Remark}

\renewcommand{\tilde}{\widetilde}
\renewcommand{\hat}{\widehat}
\newcommand{\norm}[1]{\left\| #1 \right\|}
\newcommand{\dd}{\, {\rm d}}
\newcommand{\N}{\mathbb{N}}

\newcommand{\R}{\mathbb{R}}
\newcommand{\tq}{\ |\ }
\DeclareMathOperator{\Card}{Card}
\newcommand{\boB}{\mathcal{B}}
\newcommand{\boL}{\mathcal{L}}
\newcommand{\boN}{\mathcal{N}}
\DeclareMathOperator{\Ima}{Im} \DeclareMathOperator{\sgn}{sgn}
\newcommand{\de}{{\rm d}}
\renewcommand{\geq}{\geqslant}
\renewcommand{\leq}{\leqslant}
\renewcommand{\phi}{\varphi}
\renewcommand{\epsilon}{\varepsilon}
\newcommand{\epsilo}{\epsilon_0}

\newcommand{\qedfact}{\renewcommand{\qedsymbol}{$\Diamond$}}

\DeclareMathOperator{\am}{{\alpha_{\text{min}}}}
\DeclareMathOperator{\aM}{{\alpha_{\text{max}}}}
\DeclareMathOperator{\distv}{{d_{\text{vert}}}}
\DeclareMathOperator{\dLeb}{dLeb} \DeclareMathOperator{\Leb}{Leb}
\DeclareMathOperator{\Cov}{Cov}
\newcommand{\xx}{\textbf{x}}
\newcommand{\yy}{\textbf{y}}

\title{Statistical properties of a skew product with a curve of
neutral points
  \footnote{\emph{keywords}: intermittency, countable
    Markov shift, central limit theorem, stable laws.
  \emph{2000 Mathematics Subject Classification:} 37A50, 37C40,
  60F05}}

\author{Sébastien Gouëzel
  \footnote{Département de Mathématiques et Applications,
École Normale Supérieure, 45 rue d'Ulm 75005 Paris (France).
e-mail \texttt{Sebastien.Gouezel@ens.fr}}}

\date{November 2003}

\begin{document}
\maketitle

\begin{abstract}
We study a skew product with a curve of neutral points. We show
that there exists a unique absolutely continuous invariant
probability measure, and that the Birkhoff averages of a
sufficiently smooth observable converge to a normal law or a
stable law, depending on the average of the observable along the
neutral curve.
\end{abstract}

\section{Introduction}

Let $T:M\to M$ be a map on a compact space. While uniformly
hyperbolic or uniformly expanding dynamics are well understood,
problems arise when there are neutral fixed points (where the
differential of $T$ has an eigenvalue equal to $1$). The
one-dimensional case has been thoroughly studied, particularly
when $T$ has only one neutral fixed point (see
\cite{liverani_saussol_vaienti} and references therein). The normal
form at the fixed point
dictates the asymptotics of the
dynamics, and in particular the speed of mixing, and the
convergence of Birkhoff sums to limit laws
(\cite{gouezel:stable}).

In this article, we study the same type of phenomenon, but in
higher dimension. Contrary to \cite{hu:almost_hyperbolic},
\cite{pollicott_yuri:indifferent} (where the case of isolated
fixed points is considered), our models admit a whole invariant
neutral curve.
We show that the one-dimensional results
remain essentially true.

More precisely, define a map $T_\alpha$ on $[0,1]$ by
  \begin{equation*}
  T_\alpha(x)=\left\{ \begin{array}{cl}
   x(1+2^\alpha x^\alpha) &\text{if }0\leq x\leq 1/2
   \\
   2x-1 &\text{if }1/2<x\leq 1
   \end{array}\right.
  \end{equation*}
It has a neutral fixed point at $0$, behaving like
$x(1+x^\alpha)$. To mix different such behaviors, we consider a
skew product, similar to the Alves-Viana map
(\cite{viana:multidim_attr}) but where the unimodal maps are
replaced by $T_\alpha$. Let $\alpha:S^1 \to (0,1)$ be a map with
minimum $\am$ and maximum $\aM$. Assume that
  \begin{enumerate}
  \item  $\alpha$ is $C^2$.
  \item  $0<\am<\aM<1$.
  \item $\alpha$ takes the value $\am$ at a unique point $x_0$, with
  $\alpha''(x_0)>0$.
  \item $\aM < \frac{3}{2}\am$ (which implies $\aM<\am+1/2$).
  \end{enumerate}
These conditions are for example satisfied by $\alpha(\omega)=\am
+\epsilon(1+\sin(2\pi \omega))$ where $\am\in (0,1)$ and
$\epsilon$ is small enough.

We define a map $T$ on $S^1 \times [0,1]$ by
  \begin{equation}
  \label{definit_T}
  T(\omega,x)= (F(\omega), T_{\alpha(\omega)}(x))
  \end{equation}
where $F(\omega)=4\omega$.

In the following, we will generalize to this skew product the
one-dimensional results on the maps $T_\alpha$. First of all, in
Section \ref{section_invariante}, we prove that there exists a
unique absolutely continuous invariant probability measure $m$,
whose density $h$ is in fact Lipschitz on every compact subset of
$S^1\times (0,1]$ (Theorem \ref{existence_mesure_invariante}). In
Section \ref{limit_Markov}, we prove limit theorems for abstract
Markov maps (using a method essentially due to
\cite{melbourne_torok} and recalled in Appendix
\ref{appendice:loi_stable}, and estimates of
\cite{aaronson_denker} and \cite{gouezel:stable}). Finally, in
Sections \ref{section_estimee_Xn} and \ref{section:limite}, we
study the limit laws of Birkhoff sums for the skew product $T$,
and we obtain the convergence to a normal law or a stable law,
depending on the value of $\am$. We obtain the following theorem
(see Theorem \ref{enonce_theoreme_limite} for more details).
\begin{thm*}
Set
  \begin{equation*}
  A=\frac{1}{4\left( \am^{3/2}
  \sqrt{\frac{\pi}{2\alpha''(x_0)}}\right)^{1/\am}}\int_{S^1\times \{1/2\}}
  h\dLeb,
  \end{equation*}
where $h$ is the density of the absolutely continuous invariant
probability measure.

Let $f$ be a Lipschitz function on $S^1\times [0,1]$, with $\int
f\dd m=0$. Write $c=\int_{S^1\times\{0\}} f \dLeb$ and $S_n
f=\sum_{k=0}^{n-1} f\circ T^k$. Then
  \begin{itemize}
  \item If $\am<1/2$, there exists $\sigma^2 \geq 0$ such that
$\frac{1}{\sqrt{n}} S_n f \to \boN(0,\sigma^2)$.
  \item If $\am=1/2$ and $c\not=0$, then
$\frac{S_n f}{\sqrt{ \frac{c^2 A }{4} n (\ln n)^2}} \to
\boN(0,1)$.
  \item If $1/2<\am<1$ and $c\not =0$, then
$\frac{S_n f}{n^{\am}\sqrt{\am \ln n}} \to Z$, where the random
variable $Z$ has an explicit stable distribution.
  \item If $1/2\leq \am<1$ and $c=0$, then there exists
$\sigma^2 \geq 0$ such that $\frac{1}{\sqrt{n}} S_n f \to
\boN(0,\sigma^2)$.
  \end{itemize}
\end{thm*}
An interesting feature of this example is that its study involves
sophisticated mixing properties of $F$, particularly a multiple
decorrelation property, proved in Appendix \ref{appendice:pene}
using \cite{pene:averaging}.

Theorems of \cite{gouezel:stable} could be used instead of the
method of \cite{melbourne_torok} to get the limit laws. However,
this elementary method is interesting in its own right, and can be
generalized more easily to other settings than the results of
\cite{gouezel:stable} (in particular to the case of more neutral
fixed points).

\begin{rmq}
The specific form of $F$ is of no importance at all, the results
remain true when $F$ is $C^2$ with $|F'| \geq 4$ (for example
$F(\omega)=d\omega$ with $d \geq 4$). In the same way, the only
important properties of the maps $T_\alpha$ are their normal form
close to $0$ and the fact that they are Markov. Finally, the
hypothesis $\alpha''(x_0)\not=0$ is only useful for limit
theorems, and can be replaced by: $\exists m,
\alpha^{(m)}(x_0)\not=0$ (but the normalizing factors have to be
modified accordingly). For the sake of simplicity, we will
restrict ourselves in what follows to the aforementioned case.
\end{rmq}

In this article, $a(n) \sim b(n)$ means that $a(n)/b(n) \to 1$
when $n\to \infty$. The integral with respect to a probability
measure will sometimes be denoted by $E(\cdot)$. Finally, $\lfloor
x \rfloor$ will denote the integer part of $x$.

\section{Invariant measure}
\label{section_invariante}

An important property of the map $T$, that will be used thoroughly
in what follows, is that it is Markov: there exists a partition of
the space such that every element of this partition is mapped by
$T$ on a union of elements of this partition. In fact, we will
consider $T_Y$ (the induced map on $Y=S^1 \times (1/2,1]$), which
is also Markov, and expanding, contrary to $T$. We will apply to
$T_Y$ classical results on expanding Markov maps (also called
\emph{Gibbs-Markov} maps), which we recall in the next paragraph.

\subsection{Markov maps and invariant measures}

Let $(Y,\boB,m_Y)$ be a standard probability space, endowed with a
bounded metric $d$. A non-singular map $T_Y$ defined on $Y$ is
said to be a \emph{Markov map} if there exists a finite or
countable partition $\alpha$ of $Y$ such that $\forall a\in
\alpha$, $m_Y(a)>0$, $T_Y(a)$ is a union (mod $0$) of sets of
$\alpha$, and $T_Y:a \to T_Y(a)$ is invertible. In this case,
$\alpha$ is a \emph{Markov partition} for $T_Y$.

A Markov map $T_Y$ (with a Markov partition $\alpha$) is a
\emph{Gibbs-Markov map} (\cite{aaronson:book}) if
\begin{enumerate}
\item $T_Y$ has the big image property: $\inf_{a\in
\alpha} m_Y(T_Y(a))>0$.
\item
\label{enumere_expansion}
There exists $\lambda>1$ such that $\forall a\in \alpha, \forall
x,y\in a, d(T_Yx, T_Y y)\geq \lambda d(x,y)$.
\item
\label{enumere_distortion}
Let $g$ be the inverse of the jacobian of $T_Y$, i.e.\ on a set
$a\in \alpha$, $g(x)=\frac{\dd m_Y}{\dd \left(m_Y \circ
(T_Y)_{|a}\right)} (x)$. Then there exists $C>0$ such that for all
$a\in \alpha$, for almost all $x,y\in a$,
  \begin{equation*}
  \left|1-\frac{g(x)}{g(y)} \right| \leq C d(T_Yx,T_Yy).
  \end{equation*}
\end{enumerate}

This definition is slightly more general than the definition of
\cite{aaronson:book}: the distance $d=d_\tau$ considered there is
given by $d_\tau(x,y)=\tau^{s(x,y)}$ where $\tau<1$ and $s(x,y)$
is the separation time of $x$ and $y$, i.e.\   \begin{equation}
  \label{definit_separation}
  s(x,y)=\inf\{ n\in \N \tq \nexists a\in \alpha, T^n x\in a, T^n y\in
  a\}.
  \end{equation}

The proof of \cite[Theorem 4.7.4]{aaronson:book} still works in
our context, and gives:

\begin{thm}
\label{existe_mesure_invariante}
Let $T_Y$ be a transitive Gibbs-Markov map ($\forall a,b \in
\alpha, \exists n\in \N, m_Y(T_Y^n a \cap b)>0$) such that
$\Card(\alpha_*)<\infty$, where $\alpha_*$ is the partition
generated by the images $T_Y(a)$ for $a\in \alpha$. Then $T_Y$ is
ergodic, and there exists a unique absolutely continuous (with
respect to $m_Y$) invariant probability measure, denoted by
$\mu_Y$.

Moreover, $\mu_Y=h m_Y$ where the density $h$ is bounded and
bounded away from $0$, and Lipschitz on every set of $\alpha_*$.
\end{thm}

\subsection{Preliminary estimates}

To apply Theorem \ref{existe_mesure_invariante}, we will construct
a Markov partition, and control the distortion of the inverse
branches of $T_Y$.

We will write $T_\omega^n=T_{\alpha(F^{n-1}\omega)}
\circ\cdots\circ T_{\alpha(\omega)}$, whence
$T^n(\omega,x)=(F^n\omega, T_\omega^n(x))$. Write also
$d((\omega_1,x_1),(\omega_2,x_2))= |\omega_1-\omega_2|+|x_1-x_2|$.
A point of $S^1 \times [0,1]$ will be denoted by $\xx=(\omega,x)$.
Finally, set $\distv((\omega_1,x_1),(\omega_2,x_2))=|x_2-x_1|$.

Define $X_0(\omega)=1$, $X_1(\omega)=1/2$, and for $n\geq 2$,
$X_n(\omega)$ is the preimage in $[0,1/2]$ of $X_{n-1}(F\omega)$
by $T_{\alpha(\omega)}$. These $X_n$ will be useful in the
construction of a Markov partition for $T$ (paragraph
\ref{construit_markov}).

\begin{prop}
\label{estime_croissance_grossiere_Xn}
There exists $C>0$ such that $\forall n\in \N^*, \forall \omega\in
S^1$,
  \begin{equation*}
  \frac{1}{Cn^{1/\am}} \leq X_n(\omega) \leq \frac{C}{n^{1/\aM}}.
  \end{equation*}
\end{prop}
\begin{proof}
Write $Z_1=1/2$ and $T(Z_{n+1})=Z_n$ where $T(x)=x(1+2^{\aM}
x^{\am})$. We easily check inductively that
 $Z_n\leq X_n(\omega)$ for every
$\omega$, since $T(x)\geq T_{\alpha(\omega)}(x)$ for every
$\omega$. It is thus sufficient to estimate $Z_n$ to get the
minoration. As $T(x)\geq x$, the sequence $Z_n$ is decreasing, and
nonnegative, whence it tends to a fixed point of $T$, necessarily
$0$.

We have
  \begin{align*}
  \frac{1}{Z_n^{\am}}&
  =\frac{1}{Z_{n+1}^{\am}}\left(1+2^{\aM} Z_{n+1}^{\am}\right)^{-\am}
  =\frac{1}{Z_{n+1}^{\am}}\left(1-\am 2^{\aM} Z_{n+1}^{\am}
  +o(Z_{n+1}^{\am})\right)
  \\&
  =\frac{1}{Z_{n+1}^{\am}}-\am 2^{\aM}+o(1).
  \end{align*}
A summation gives $\frac{1}{Z_m^{\am}}\sim m \am 2^{\aM}$, whence
$Z_m \sim C/m^{1/\am}$, which concludes the minoration.

The majoration is similar, using a sequence $Z'_n$ with $Z'_n \geq
X_n(\omega)$.
\end{proof}


\emph{We fix once and for all a large enough constant $D$.} The
following definition is analogous to a definition of Viana
(\cite{viana:multidim_attr}).

\begin{defn}
Let $\psi: K\to [0,1]$, where $K$ is a subinterval of $S^1$. We
say that the graph of $\psi$ is an admissible curve if $\psi$ is
$C^1$ with $|\psi'|\leq D$.
%
\end{defn}

\begin{prop}
Let $\psi$ be an admissible curve, defined on $K$ with $|K| <1/4$,
and included in $K\times [0,1/2]$ or $K\times (1/2,1]$. Then the
image of $\psi$ by $T$ is still an admissible curve.
\end{prop}
\begin{proof}
Let $(u,v)$ be a tangent vector at $(\omega,x)$ with $|v|\leq D
|u|$, we have to check that its image $(u',v')$ by $DT(\omega,x)$
still satisfies $|v'|\leq D |u'|$.

Assume first that $x\leq 1/2$, whence $u'=4u$ and
$v'=(1+(2x)^{\alpha(\omega)} (\alpha(\omega)+1))v+ x
\ln(2x)\alpha'(\omega)(2x)^{\alpha(\omega)}u$. As
$\alpha(\omega)\leq \aM\leq 1$, we get $|v'|\leq 3|v|+C|u|$ for a
constant $C$ (which depends only on $\norm{\alpha'}_\infty$).
Thus,
  \begin{equation*}
  \frac{|v'|}{|u'|} \leq  \frac{3}{4} \frac{|v|}{|u|} +\frac{C}{4}.
  \end{equation*}
This will give $|v'|/|u'| \leq D$ if $\frac{3}{4}D +\frac{C}{4}
\leq D$, which is true if $D$ is large enough.

Assume then that $x>1/2$. Then $u'=4u$ and $v'=2v$, and there is
nothing to prove.
\end{proof}

\begin{cor}
\label{controle_rapport_horz_vert}
Let $(\omega_1,x_1)$ and $(\omega_2,x_2)$ be two points in
$S^1\times [0,1/2]$ with $|x_1-x_2|\leq D|\omega_1-\omega_2|$ and
$|\omega_1-\omega_1|\leq \frac{1}{8}$. Then their images satisfy
$|x'_1-x'_2|\leq D|\omega'_1-\omega'_2|$.
\end{cor}
\begin{proof}
Use a segment between the two points: it is an admissible curve,
whence its image is still admissible.
\end{proof}

\subsection{The Markov partition}

\label{construit_markov}

Set $Y=S^1 \times (1/2,1]$. For $\xx \in Y$, set
$\phi_Y(\xx)=\inf\{n>0 \tq T^n(\xx)\in Y\}$: this is the first
return time to $Y$, everywhere finite. The map
$T_Y(\xx):=T^{\phi_Y(\xx)}(\xx)$ is the map induced by $T$ on $Y$.
We will show that $T_Y$ is a Gibbs-Markov map, by constructing an
appropriate Markov partition.

If $I$ is an interval of $S^1$, we will abusively write $I\times
[X_{n+1},X_n]$ for $\{(\omega,x)\tq \omega \in I, x\in
[X_{n+1}(\omega), X_n(\omega)]\}$.

Set $I_n(\omega)=[X_{n+1}(\omega),X_n(\omega)]$ (or $\{\omega\}
\times [X_{n+1}(\omega),X_n(\omega)]$, depending on the context).
By definition of $X_n$, $T$ maps $\{\omega\}\times I_n(\omega)$
bijectively on $\{F\omega\} \times I_{n-1}(F\omega)$. Thus, the
interval $I_n(\omega)$ returns to $[1/2,1]$ in exactly $n$ steps.

Let $Y_n(\omega)$ be the preimage in $[1/2,1]$ of
$X_{n-1}(F\omega)$ under $T_{\alpha(\omega)}$. Thus, the interval
$J_n(\omega)= [Y_{n+1}(\omega),Y_n(\omega)]$ returns to $[1/2,1]$
in $n$ steps.

\emph{We fix once and for all $0<\epsilo<\frac{1}{8}$, small
enough so that $D\epsilo$ is less than the length of every
interval $I_1(\omega)$}. (This condition will be useful in
distortion estimates).

\emph{Let $q$ be large enough so that $\frac{1}{4^q}<\epsilo$, and
consider $A_{s,n}=\left[
\frac{s}{4^{q+n}},\frac{s+1}{4^{q+n}}\right] \times J_n$, for
$n\in \N^*$ and $0\leq s\leq 4^{q+n}-1$}: this set is mapped by
$T^n$ on $\left[\frac{s}{4^q},\frac{s+1}{4^q}\right] \times
[1/2,1]$. Let $K_0,\ldots,K_{4^q-1}$ be the sets
$\left[\frac{i}{4^q},\frac{i+1}{4^q}\right] \times [1/2,1]$. Then
the map $T_Y$ is an isomorphism between each $A_{s,n}$ and some
$K_i$. Consequently, the map $T_Y$ is Markov for the partition
$\{A_{s,n}\}$, and it has the big image property.

To apply Theorem \ref{existe_mesure_invariante}, we need expansion
(for \eqref{enumere_expansion} in the definition of Gibbs-Markov
maps) and distortion control (for \eqref{enumere_distortion}). The
expansion is given by the next proposition, and the distortion is
estimated in the following paragraph.

On the intervals $[X_3(\omega),X_1(\omega)]$, the derivative of
$T_{\alpha(\omega)}$ is greater than $1$, whence greater than a
constant $2>\lambda>1$, independent of $\omega$.

For $(\omega_1,x_1)$ and $(\omega_2,x_2)\in S^1\times [0,1]$, set
  \begin{equation}
  \label{definit_d'}
  d'((\omega_1,x_1),(\omega_2,x_2))=a|x_1-x_2|+|\omega_1-\omega_2|
  \end{equation}
where $a=\frac{1-\lambda/4}{D}$.

\begin{prop}
\label{dilate_markov}
On each $A_{s,n}$, the map $T^n$ is expanding by at least
$\lambda$ for the distance $d'$.
\end{prop}
\begin{proof}
For $n=1$ (the points return directly to $S^1\times [1/2,1]$),
everything is linear, and the result is clear. Assume $n\geq 2$.

Take $(\omega_1,x_1)$ and $(\omega_2,x_2)\in A_{s,n}$, with for
example $x_2\geq x_1$. The points $(\omega_1,x_1)$ and
$(\omega_2,x_1)$ return to $S^1\times [1/2,1]$ after at least $n$
iterations (by hypothesis for the first point, and the second
point is under $(\omega_2,x_2)$). We can use Corollary
\ref{controle_rapport_horz_vert} $n-1$ times, and get that in
vertical distance, $\distv( T^n(\omega_1,x_1),T^n(\omega_2,x_1))
\leq D|F^n\omega_1-F^n\omega_2|$. In particular,
$T_{\omega_2}^n(x_1) \geq T_{\omega_1}^n(x_1) - D \epsilon_0 \geq
1/2-D \epsilon_0$. Thus, by definition of $\epsilon_0$,
$T^n(\omega_2,x_1) \in I_i(F^n \omega_2)$ for $i=0$ or $1$, whence
$T^{n-1}(\omega_2,x_1) \in [X_3(F^{n-1} \omega_2),
X_1(F^{n-1}\omega_2)]$. Note that $T^{n-1}(\omega_2,x_2)$ belongs
to the same interval (in fact, $T^{n-1}_{\omega_2}(x_2) \in
[X_2(F^{n-1}\omega_2), X_1(F^{n-1}\omega_2)]$). Moreover, the
$T_\alpha$ are expanding, whence $\distv(T^{n-1}(\omega_2,x_1),
T^{n-1}(\omega_2,x_2)) \geq |x_1-x_2|$. We apply once more $T$,
which expands at least by $\lambda$ on $[X_3(F^{n-1}\omega_2),
X_1(F^{n-1} \omega_2)]$ by definition of $\lambda$, and get
$\distv(T^n(\omega_2,x_1),T^n(\omega_2,x_2)) \geq
\lambda|x_1-x_2|$.

Finally,
  \begin{align*}
  d'(T^n(\omega_1,x_1),T^n(\omega_2,x_2))&
  =a
  \distv(T^n(\omega_1,x_1),T^n(\omega_2,x_2))+|F^n\omega_1-F^n\omega_2|
  \\&
  \geq a \distv(T^n(\omega_2,x_1),T^n(\omega_2,x_2))-a \distv(
  T^n(\omega_1,x_1),T^n(\omega_2,x_1))
  \\& \hphantom{=\ }+|F^n\omega_1-F^n\omega_2|
  \\&
  \geq a \lambda|x_1-x_2| -aD
  |F^n\omega_1-F^n\omega_2|+|F^n\omega_1-F^n\omega_2|.
  \end{align*}
The proposition will be proved if $(1-aD)|F^n\omega_1-F^n\omega_2|
\geq \lambda |\omega_1-\omega_2|$. But
  \begin{equation*}
  (1-aD)|F^n\omega_1-F^n\omega_2|
  =(1-aD)4^n |\omega_1-\omega_2|
  \geq (1-aD)4 |\omega_1-\omega_2|
  =\lambda |\omega_1-\omega_2|.
  \end{equation*}
\end{proof}

\subsection{Distortion bounds}

\begin{lem}
\label{lemme_distortion_deplac_horz}
There exists a constant $E>0$ such that
  \begin{equation}
  \begin{split}
  \forall n>0, \forall \omega_1,\omega_2&\in S^1 \text{ with
  }|\omega_1-\omega_2|\leq \frac{\epsilo}{4^n}, \forall x_1 \in
  J_n(\omega_1) \text{ with }T_{\omega_2}^{n-1} x_1\leq 1/2,\\&
  \left| \ln (T_{\omega_1}^n)'(x_1)-\ln (T_{\omega_2}^n)'(x_1)\right|
  \leq E |F^n \omega_1 -F^n \omega_2|.
  \end{split}
  \end{equation}
\end{lem}
\begin{proof}
We use Corollary \ref{controle_rapport_horz_vert} $n-1$ times  and
get for $0\leq k\leq n$ that $|T_{\omega_1}^k x_1 -T_{\omega_2}^k
x_1|\leq D|F^k\omega_1-F^k\omega_2|$.

In particular, for $k=n$, $|T_{\omega_1}^n x_1|\geq 1/2$, whence
$|T_{\omega_2}^n x_1|\geq 1/2-D\epsilo$. Consequently,
$T^n(\omega_2,x_1)\in I_i(F^n\omega_2)$ for some $i\in \{0,1\}$,
by definition of $\epsilo$. An inverse induction gives
$T^k(\omega_2,x_1) \in I_{n-k+i}(F^k\omega_2)$.

For $x\leq 1/2$ and $\omega\in S^1$, write $G(\omega,x)=\ln
T_{\alpha(\omega)}'(x)=
\ln\left(1+(\alpha(\omega)+1)(2x)^{\alpha(\omega)}\right)$. Then
  \begin{equation*}
  \frac{\partial G}{\partial x}(\omega,x)
  =\frac{(\alpha(\omega)+1)\alpha(\omega)2^{\alpha(\omega)}
  x^{\alpha(\omega)-1}}{1+(\alpha(\omega)+1)(2x)^{\alpha(\omega)}}
  \leq C x^{\am-1}
  \end{equation*}
and
  \begin{equation*}
  \left|\frac{\partial G}{\partial \omega}(\omega,x)\right|
  =\left|\frac{\alpha'(\omega)(2x)^{\alpha(\omega)}
  +(\alpha(\omega)+1)\alpha'(\omega)\ln(2x)(2x)^{\alpha(\omega)}}
  {1+(\alpha(\omega)+1)(2x)^{\alpha(\omega)}} \right|
  \leq C.
  \end{equation*}

Lemma \ref{estime_croissance_grossiere_Xn}, and the fact that
$T^k(\omega_1,x_1)\in I_{n-k}(F^k\omega_1)$ and $T^k(\omega_2,x_1)
\in I_{n-k+i}(F^k\omega_2)$ with $i\leq 1$, give that the second
coordinates of $T^k(\omega_1,x_1)$ and $T^k(\omega_2,x_1)$ are
$\geq \frac{1}{C (n-k+1)^{1/\am}}$. On the set of points
$(\omega,x)$ with $x\geq \frac{1}{C (n-k+1)^{1/\am}}$, the
estimates on the partial derivatives of $G$ show that this
function is $C (n-k+1)^{1/\am -1}$-Lipschitz, whence
  \begin{align*}
  |G(T^k(\omega_1,x_1))-G(T^k(\omega_2,x_1))| &
  \leq C (n-k+1)^{1/\am -1} d((T^k(\omega_1,x_1),T^k(\omega_2,x_1))
  \\&
  \leq C (n-k+1)^{1/\am -1} (1+D)|F^k \omega_1-F^k \omega_2|
  \\&
  \leq C (n-k+1)^{1/\am -1} (1+D) 4^k |\omega_1-\omega_2|.
  \end{align*}
Finally,
  \begin{align*}
  \left| \ln (T_{\omega_1}^n)'(x_1)-\ln (T_{\omega_2}^n)'(x_1)\right|&
  \leq \sum_{k=0}^{n-1}|G(T^k(\omega_1,x_1))-G(T^k(\omega_2,x_1))|
  \\&
  \leq C 4^n |\omega_1-\omega_2| \sum_{k=0}^{n-1} (n-k+1)^{1/\am
  -1} 4^{k-n}
  \\&
  \leq C |F^n\omega_1-F^n\omega_2| \sum_{l=1}^\infty (l+1)^{1/\am
  -1} 4^{-l}.
  \end{align*}
The last sum is finite, which concludes the proof.
\end{proof}

For $n\geq 2$, write
$J_n^+(\omega)=[Y_{n+2}(\omega),Y_n(\omega)]$. Thus, if $n\geq 1$,
$J_{n+1}^+(\omega)$ is the preimage of $I_n^+(F\omega)$, defined
by $I_n^+(F\omega)=[X_{n+2}(F\omega),X_{n}(F\omega)]$. These
intervals will appear naturally in distortion controls, since we
have seen in the proof of Lemma \ref{lemme_distortion_deplac_horz}
that, if we move away horizontally from a point of
$J_n(\omega_1)$, we find a point of $J_{n+i}(\omega_2)$ for $i\in
\{0,1\}$, i.e.\ in $J_n^+(\omega_2)$.

%
%

\begin{lem}
\label{distortion_verticale_bornee}
There exists a constant $C$ such that
  \begin{equation*}
  \forall n\geq 0, \forall \omega\in S^1, \forall x,y\in
  J_n^+(\omega),\ \left|\ln (T_\omega^n)'(x)-\ln(T_\omega^n)'(y)
  \right| \leq C |T_\omega^n(x)-T_\omega^n(y)|.
  \end{equation*}
\end{lem}
\begin{proof}
Recall that the Schwarzian derivative of an increasing
diffeomorphism $g$ of class $C^3$ is
$Sg(x)=\frac{g'''(x)}{g'(x)}-\frac{3}{2}\left(\frac{g''(x)}{g'(x)}
\right)^2$. The composition of two functions with nonpositive
Schwarzian derivative still has a nonpositive Schwarzian
derivative.

For $\tau>0$, the Koebe principle (\cite[Theorem
IV.1.2]{demelo_vanstrien}) states that, if $Sg\leq 0$, and
$J\subset J'$ are two intervals such that $g(J')$ contains a
$\tau$-scaled neighborhood of $g(J)$ (i.e.\ the intervals on the
left and on the right of $g(J)$ in $g(J')$ have length at least
$\tau |g(J)|$), then there exists a constant $K(\tau)$ such that
  \begin{equation*}
  \forall x,y\in J, \left|\ln g'(x) - \ln g'(y)\right| \leq
  K(\tau) \frac{|x-y|}{|J|}.
  \end{equation*}
This implies that the distortion of $g$ is bounded on $J$, whence
it is possible to replace the bound on the right with $K'(\tau)
\frac{|g(x)-g(y)|}{|g(J)|}$.

In our case, if $0<\alpha<1$, the left branch of $T_\alpha$ has
nonpositive Schwarzian derivative, since $T_\alpha'''<0$ and
$T_\alpha'>0$. Let in particular $g$ be the composition of the
left branches of $T_{\alpha(F^{n-1}\omega)},\ldots,
T_{\alpha(F\omega)}$, and of the right branch of
$T_{\alpha(\omega)}$. Then, on $J_n^+$, we have $T_\omega^n=g$,
and $g$ has nonpositive Schwarzian derivative.

We want to see that $\left|\ln (T_\omega^n)'(x) - \ln
(T_\omega^n)'(y)\right|\leq C|T_\omega^n(x)-T_\omega^n(y)|$. For
this, we apply the Koebe principle to $J=J_n^+$ and
$J'=[1/2+\delta,2]$ for $\delta$ very small. Then $g(J)=[X_2,1]$
while $g(J')$ contains $[\delta',2]$ for $\delta'>0$, arbitrarily
small if $\delta$ is small enough. As the $X_{2}$ are uniformly
bounded away from $0$, there exists $\tau>0$ (independent of
$\omega$ and $n$) such that $g(J')$ contains a $\tau$-scaled
neighborhood of $g(J)$. The Koebe principle then gives the desired
result.
\end{proof}

\begin{prop}
\label{prop_distortion_bornee}
There exists a constant $C$ such that, for every $A_{s,n}$, for
every $(\omega_1,x_1)$ and $(\omega_2,x_2)\in A_{s,n}$,
  \begin{equation*}
  \left| \frac{\det DT^n(\omega_1,x_1)}{\det DT^n(\omega_2,x_2)}-1
  \right|\leq C d(T^n(\omega_1,x_1),T^n(\omega_2,x_2)).
  \end{equation*}
\end{prop}
\begin{proof}
The matrix $DT^n(\omega,x)$ is upper triangular, with $4^n$ in the
upper left corner. Thus, we have to show that
  \begin{equation*}
  \left| \ln(T_{\omega_1}^n)'(x_1)-\ln(T^n_{\omega_2})'(x_2)
  \right|\leq C d(T^n(\omega_1,x_1),T^n(\omega_2,x_2)).
  \end{equation*}
Assume for example $x_2\geq x_1$, which implies that
$T_{\omega_2}^k(x_1)\leq 1/2$ for $k=0,\ldots,n-1$. Lemma
\ref{lemme_distortion_deplac_horz} can be applied to $x_1$,
$\omega_1$ and $\omega_2$, and gives in particular that $x_1\in
J_n^+(\omega_2)$.

Write
  \begin{align*}
  \left| \ln(T_{\omega_2}^n)'(x_2)-\ln(T^n_{\omega_1})'(x_1)
  \right|
  &\leq
  \left| \ln(T_{\omega_2}^n)'(x_2)-\ln(T^n_{\omega_2})'(x_1)
  \right|
  +  \left| \ln(T_{\omega_2}^n)'(x_1)-\ln(T^n_{\omega_1})'(x_1)
  \right|
  \\&
  \leq C
  d(T^n(\omega_2,x_2)),T^n(\omega_2,x_1))
  + E|F^n\omega_2-F^n \omega_1|
  \end{align*}
by Lemmas \ref{lemme_distortion_deplac_horz} and
\ref{distortion_verticale_bornee}. For the first term,
  \begin{align*}
  d(T^n(\omega_2,x_2),T^n(\omega_2,x_1))
  &\leq
  d(T^n(\omega_2,x_2),T^n(\omega_1,x_1))
  +d(T^n(\omega_1,x_1),T^n(\omega_2,x_1))
  \\&
  \leq d(T^n(\omega_2,x_2),T^n(\omega_1,x_1))
  + (D+1) |F^n \omega_1-F^n \omega_2|
  \end{align*}
using admissible curves.

As $|F^n \omega_1-F^n \omega_2|\leq
d(T^n(\omega_1,x_1),T^n(\omega_2,x_2))$, we get the conclusion.
\end{proof}

\subsection{Construction of the invariant measure}

The previous estimates and Theorem \ref{existe_mesure_invariante}
easily give that $T_Y$ admits an invariant measure, with Lipschitz
density. Inducing gives an invariant measure for $T$, whose
density is Lipschitz on each set $S^1\times (X_{n+1},X_n)$.
However, this does not exclude discontinuities on $S^1\times X_n$,
which is not surprising since $T$ itself has a discontinuity on
$S^1 \times \{1/2\}$, which will then propagate to the other
$X_n$, since the measure is invariant.

However, in the one-dimensional case, Liverani, Saussol and
Vaienti (\cite{liverani_saussol_vaienti}) have proved that the
density is really continuous everywhere, since they constructed it
as an element of a cone of continuous functions. This fact remains
true here:

\begin{thm}
\label{existence_mesure_invariante}
The map $T$ admits a unique absolutely continuous invariant
probability measure $\dd m$. Moreover, this measure is ergodic.
Finally, the density $h=\frac{\dd m}{\dLeb}$ is Lipschitz on every
compact subset of $S^1 \times (0,1]$.
\end{thm}

\begin{proof}
Consider the map $T_Y$ induced by $T$ on $Y=S^1\times (1/2,1]$. It
is Markov for the partition $\alpha=\{A_{s,n}\}$, and transitive
for this partition since $T_Y^2(a)=Y$ for all $a\in \alpha$.
Moreover, it is expanding for $d'$ on each set of the partition
(Proposition \ref{dilate_markov}) and its distortion is Lipschitz
(Proposition \ref{prop_distortion_bornee}, and $d$ equivalent to
$d'$).

Theorem \ref{existe_mesure_invariante} shows that $T_Y$ admits a
unique absolutely continuous invariant probability measure $\dd
m_Y=h \dLeb$, which is ergodic. Moreover, the density $h$ is
Lipschitz (for the distance $d'$, whence for the usual one) on
each element of the partition $\alpha_*$ generated by the sets
$T_Y(a)$, i.e.\ on the sets $K_i$.

To construct an invariant measure for the initial map $T$, we use
the classical induction process (\cite[Section
1.1.5]{aaronson:book}): let $\phi_Y$ be the return time to $Y$
under $T$, then $\mu=\sum_{n=0}^\infty T_*^n(m_Y | \phi_Y>n)$ is
invariant. To check that the new measure has finite mass, we have
to see that $\sum m_Y(\phi_Y>n) <\infty$. As $\dd m_Y$ and $\dLeb$
are equivalent, we check it for $\dLeb$. We have
  \begin{equation*}
  \Leb(\phi_Y>n)=\Leb(S^1 \times [1/2,Y_{n+1}])=\frac{1}{2}
  \Leb(S^1\times [0,X_n])
  \leq \frac{1}{2} \frac{C}{n^{1/\aM}},
  \end{equation*}
using Lemma \ref{estime_croissance_grossiere_Xn}. As $\aM<1$, this
is summable.

We know that $h$ is Lipschitz on the sets
$[\frac{s}{4^q},\frac{s+1}{4^q}]\times [1/2,1]$, we have to prove
the continuity on $\{s/4^q\}\times [1/2,1]$, which is not hard:
these numbers $s/4^q$ are artificial, since they depend on the
arbitrary choice of a Markov partition on $S^1$. We can do the
same construction using other sets than the $A_{s,n}$. For
example, set $A'_{s,n}=\left[ \frac{1}{3}+\frac{s}{4^{q+n}},
\frac{1}{3}+\frac{s+1}{4^{q+n}}\right]\times J_n$, and
$K'_i=\left[\frac{1}{3}+\frac{i}{4^q},
\frac{1}{3}+\frac{i+1}{4^q}\right]$. Since $1/3$ is a fixed point
of $F$, the map $T_Y$ is Markov for the partition $\{A'_{s,n}\}$,
and each of these sets is mapped on a set $K'_i$. Thus, the same
arguments as above apply, and prove that $h$ is Lipschitz on each
set $K'_i$. Since the boundaries of the sets $K_i$ and $K'_i$ are
different, this shows that $h$ is in fact Lipschitz on $S^1\times
[1/2,1]$.

We show now that $h$ is Lipschitz on $S^1 \times [X_2,1]$. Note
that it is slightly incorrect to say that $h$ is Lipschitz, since
$h$ is defined only almost everywhere. Nevertheless, if we prove
that $|h(\xx)-h(\yy)|\leq Cd(\xx,\yy)$ for almost all $\xx$ and
$\yy$, then there will exist a unique version of $h$ which is
really Lipschitz. Thus, all the equalities we will write until the
end of this proof will be true only almost everywhere.

Let $A_{s,n}^+ = \left[\frac{s}{4^{q+n}},
\frac{s+1}{4^{q+n}}\right]\times J_n^+$: $T^n$ is a diffeomorphism
between $A_{s,n}^+$ and $K_i^+ =
\left[\frac{i}{4^q},\frac{i+1}{4^q}\right]\times [X_2,1]$. We fix
some $K^+=K_i^+=I\times [X_2,1]$, and we show that $h$ is
Lipschitz on $K^+$. Let $U_1,U_2,\ldots$ be the inverse branches
of $T^{n_1},T^{n_2},\ldots$ whose images all coincide with $K^+$.
Let $T_Y$ be the map induced by $T$ on $Y=S^1\times [1/2,1]$. Then
$h\dLeb_{|Y}$ is invariant under $T_Y$, which means that, for each
$\xx\in I\times [1/2,1]$,
  \begin{equation*}
  h(\xx)=\sum JU_j(\xx) h(U_j \xx)
  \end{equation*}
where $JU_j$ is the jacobian of $U_j$.

Let $Z=S^1\times [X_2,1]$, and $T_Z$ be the map induced by $T$ on
$Z$. Since $h\dLeb_{|Z}$ is also invariant under $T_Z$, we have
the same kind of equation as above. For $\xx\in I\times
[X_2,1/2]$, all its preimages under $T_Z$ are in $S^1 \times
[1/2,1]$, and the invariance gives that
  \begin{equation*}
  h(\xx)=\sum JU_j(\xx) h(U_j \xx).
  \end{equation*}

We have shown that, for every $\xx\in S^1 \times [X_2,1]$,
  \begin{equation*}
  h(\xx)=\sum JU_j(\xx) h(U_j \xx).
  \end{equation*}
This means that $h$ is invariant under some kind of transfer
operator, even though it is not a real transfer operator since the
images of the maps $U_j$ are not disjoint, and since they do not
cover the space. In particular, the images of the $U_j$ are
included in $S^1 \times [1/2,1]$, and we already know that $h$ is
Lipschitz on this set.

The bounds of the previous paragraphs still apply to the
distortion of the $U_j$, and their expansion. In particular,
$\left|1-\frac{JU_j(\yy)}{JU_j(\xx)}\right|\leq Cd(\xx,\yy)$ for a
constant $C$ independent of $j$, and $|h(U_j \xx)-h(U_j \yy)| \leq
C d(U_j \xx,U_j \yy)\leq D d(\xx,\yy)$ (since $h$ is Lipschitz on
the image of $U_j$). Thus,
  \begin{align*}
  |h(\xx)-h(\yy)|&
  \leq \sum |JU_j(\xx)h(U_j \xx)-JU_j(\yy)h(U_j \yy)|
  \\&
  \leq \sum |JU_j(\xx)| \left|1-\frac{JU_j(\yy)}{JU_j(\xx)}\right|
  |h(U_j
  \xx)| +\sum |JU_j(\yy)| |h(U_j \xx)-h(U_j \yy)|
  \\&
  \leq Cd(\xx,\yy)\sum |JU_j(\xx)| + Dd(\xx,\yy) \sum |JU_j(\yy)|.
  \end{align*}
It remains to prove that $\sum |JU_j(\xx)|$ is bounded. The bound
on distortion gives $JU_j(\xx) \asymp \Leb(\Ima U_j)$, whence
$\sum JU_j(\xx)\leq C\sum \Leb(\Ima U_j)$, which is finite since
every point of $I\times[1/2,1]$ is in the image of at most two
maps $U_j$.

We have proved that $h$ is Lipschitz on $S^1\times [X_2,1]$,
except maybe on $\{\frac{s}{4^q}\}\times [X_2,1]$. As above, using
another Markov partition, we exclude the possibility of
discontinuities there. Thus, $h$ is Lipschitz on $S^1 \times
[X_2,1]$.

%

To prove that $h$ is Lipschitz on $S^1\times [X_k,1]$, we do
exactly the same thing, except that we consider $[Y_{n+k},Y_n]$
instead of $J_n^+=[Y_{n+2},Y_n]$. As above, writing
$U_1,U_2,\ldots$ for the inverse branches of $T^n$ defined on a
set $[\frac{s}{4^{n+q}},\frac{s+1}{4^{n+q}}]\times [Y_{n+k},Y_n]$
and whose image is $K'=[\frac{i}{4^q},\frac{i+1}{4^q}] \times
[X_k,1]=I\times [X_k,1]$, we show that $h(\xx)=\sum JU_j(\xx)
h(U_j\xx)$ for $\xx\in K'$. In fact, for $\xx \in I\times
[X_{l},X_{l-1}]$, we use the invariance of $h\dLeb$ under the map
induced by $T$ on $S^1 \times [X_l,1]$. We conclude finally as
above, using the fact that $h$ is Lipschitz on $S^1\times
[1/2,1]$, which contains the images of the $U_j$.

This concludes the proof, since every compact subset of
$S^1\times(0,1]$ is contained in $S^1 \times [X_k,1]$ for large
enough $k$.
\end{proof}

\section{Limit theorems for Markov maps}
\label{limit_Markov}

We want to establish limit theorems for Birkhoff sums, of the form
$\sum_{k=0}^{n-1} f(T^k x)$. We give in this section an abstract
result, valid for a map that induces a Gibbs-Markov map on a
subset of the space (which is the case of our skew product).
Related limit theorems have been proved in \cite{gouezel:stable},
but we will show here a slightly different result, which requires
more control on the return time $\phi$ but is more elementary,
using Theorem \ref{thm_probabiliste_general} proved in Appendix
\ref{appendice:loi_stable} and inspired by results of Melbourne
and Török (\cite{melbourne_torok}) for flows. An advantage of this
new method is that, contrary to \cite{gouezel:stable}, it can
easily be extended to stable laws of index $1$.

If $Z_0,\ldots,Z_{n-1},\ldots$ are independent identically
distributed random variables with zero mean, the sums
$\frac{1}{B_n} \sum_{k=0}^{n-1} Z_k$ (where $B_n$ is a real
sequence) converge to a nontrivial limit in essentially three
cases: if $Z_k \in L^2$, there is convergence to a normal law for
$B_n=\sqrt{n}$. There is also convergence to a normal law, but with a
different normalization, if $P(|Z_k|>x)=x^{-2}l(x)$ with
$L(x):=2\int_1^x \frac{l(u)}{u}\dd u$ unbounded and slowly varying
(i.e.\ $L:(0,\infty) \to (0,\infty)$ satisfies $\forall a
>0, \lim_{x\to \infty} L(ax)/L(x)=1$) -- this is in particular true
when $l$ itself is slowly varying.
Finally,
if $P(Z_k>x)=(c_1+o(1))x^{-p}L(x)$ and
$P(Z_k<-x)=(c_2+o(1))x^{-p}L(x)$, where $L$ is slowly varying and
$p\in (0,2)$, we have convergence (for a good choice of $B_n$) to
a limit law called stable law. Moreover, these are the only cases
where there is a convergence (\cite{feller:2}).

In the dynamical setting, we will prove the same kind of limit
theorems, still with three possible cases: $L^2$, normal
nonstandard, and stable. The normalizations will moreover be the
same as in the probabilistic setting.

\begin{thm}
\label{thm_abstrait_markov}
Let $T:X\to X$ be an ergodic transformation preserving a
probability measure $m$. Assume that there exists a subset $Y$ of
$X$ with $m(Y)>0$ such that the first return map
$T_Y(x)=T^{\phi(x)}(x)$ (where $\phi(x)= \inf\{ n>0 \tq T^n(x) \in
Y\}$) is Gibbs-Markov for $m_{|Y}$, a partition $\alpha$ of $Y$ such
that $\phi$ is constant on each element of $\alpha$,
and a distance $d$ on $Y$.

Let $f:X \to \R$ be an integrable map with $\int f=0$, such that
$f_Y(y):=\sum_{n=0}^{\phi(y)-1} f(T^n y)$ satisfies
  \begin{equation}
  \label{condition_markov_abstrait}
  \sum_{a\in \alpha} m(a) D f_Y(a) <\infty
  \end{equation}
where
  \begin{equation*}
  D f_Y(a)=\inf\{ C >0 \tq \forall x,y\in a, |f_Y(x)-f_Y(y)| \leq C
d(x,y)\}.
  \end{equation*}

Set $M(y)=\max_{1\leq k \leq \phi(y)} \left| \sum_{j=0}^{k-1}
f\circ T^j(y) \right|$.

Then:
\begin{itemize}
\item Assume that $f_Y \in L^2$ and $M \in L^2$. Assume moreover
that $\phi$ satisfies one of the following hypotheses:
  \begin{itemize}
  \item $\phi\in L^2$.
  \item $m(\phi> x)=x^{-p} L(x)$ where $L$ is slowly varying and
  $p \in (1,2]$.
  \end{itemize}
Then there exists $\sigma^2 \geq 0$ such that $\frac{1}{\sqrt{n}}
S_n f \to \boN(0,\sigma^2)$.

\item
Assume that $m(|f_Y| > x)=x^{-2} l(x)$, with
$L(x):=2\int_1^x\frac{l(u)}{u}\dd u$ unbounded and slowly varying.
Assume moreover that
$m(M
> x) \leq C x^{-2} l(x)$, and $m(\phi>x)=(c+o(1)) x^{-2}l(x)$.
Let $B_n\to \infty$ satisfy $n L(B_n)=B_n^2$. Then
$\frac{1}{B_n} S_n f \to \boN(0,1)$.

\item
Assume that $m(f_Y>x)=(c_1+o(1))x^{-p}L(x)$ and
$m(f_Y<-x)=(c_2+o(1))x^{-p} L(x)$ where $L$ is a slowly varying
function, $p\in (1,2)$, and $c_1,c_2\geq 0$ with $c_1+c_2>0$.
Assume also that $m(M>x) \leq C x^{-p} L(x)$, and
$m(\phi>x)=(c_3+o(1)) x^{-p}L(x)$. Let $B_n\to \infty$ satisfy $n
L(B_n)=B_n^p$. Then $\frac{1}{B_n} S_n f \to Z$ where the random
variable $Z$ has a characteristic function given by
  \begin{equation*}
  E(e^{itZ})=e^{-c|t|^p \left( 1-i\beta \sgn(t) \tan
  \left(\frac{p\pi}{2} \right) \right)}
  \end{equation*}
with $c=(c_1+c_2) \Gamma(1-p) \cos \left(\frac{p\pi}{2} \right)$
and $\beta=\frac{c_1-c_2}{c_1+c_2}$.
\end{itemize}
\end{thm}

Note that $M(y)\leq \sum_{j=0}^{\phi(y)-1} |f(T^j y)|=|f|_Y(y)$.
Thus, if the integrability hypotheses of the theorem are satisfied
by $|f|_Y$ (which will often be the case), they are automatically
satisfied by $M$.

In the second case of the theorem, when $l$ itself is slowly varying,
then $L$ is automatically slowly varying.

The second case of the theorem is not the most general possible
result, since one may have convergence to a normal law even when the
function $l$ is not slowly varying (what really matters is that $L$ is 
slowly varying). The theorem can be extended without problem to this
more general setting, but the result becomes more complicated to state. In
the applications, the statement given in Theorem
\ref{thm_abstrait_markov} will be sufficient.

\begin{proof}
The idea is to use Theorem \ref{thm_probabiliste_general}: we have
to check all its hypotheses. We will use the notations of this
theorem, and in particular write $E_Y(u)=\frac{\int_Y u \dd
m}{m(Y)}$.

We first treat the third case (stable law), using the results of
\cite{aaronson_denker} (and the generalizations of
\cite{gouezel:stable}). Let $s(x,y)$ be the separation time of $x$
and $y$ defined in \eqref{definit_separation}, $\tau=1/\lambda$
and $d_\tau=\tau^s$ the corresponding metric. Since every
iteration of $T_Y$ expands by at least $\lambda$, we get
$d(x,y)\leq C d_\tau(x,y)$. In particular, we can assume without
loss of generality that $d=d_\tau$, which is the setting of
\cite{aaronson_denker} and \cite{gouezel:stable}.

Let $P$ be the transfer operator associated to $T_Y$ (i.e.\
defined by $\int u\cdot v\circ T_Y = \int P(u)\cdot v$), and
$P_t(u)=P(e^{itf_Y} u)$. Let $\boL$ be the space of bounded
Lipschitz functions (i.e.\ such that there exists $C$ such that,
$\forall a\in \alpha, \forall x,y\in a, |g(x)-g(y)|\leq C
d(x,y)$). Theorem 5.1 of \cite{aaronson_denker} ensures that, for
small enough $t$, $P_t$ acting on $\boL$ has an eigenvalue
$\lambda(t)=e^{-\frac{c}{m(Y)}|t|^p \left( 1-i\beta \sgn(t) \tan
\left(\frac{p\pi}{2} \right) \right)L(|t|^{-1})(1+o(1))}$, the
remaining part of its spectrum being contained in a disk of radius
$\leq 1-\delta<1$. In fact, this theorem requires that $Df_Y(a)$
is bounded, but \cite[Theorem 3.8]{gouezel:stable} shows that it
remains true under the weaker assumption $\sum m(a)D
f_Y(a)<\infty$.

The slow variation of $L$ easily implies that
$\lambda\left(\frac{t}{B_n} \right)^{\lfloor n m(Y) \rfloor} \to
e^{-c|t|^p \left( 1-i\beta \sgn(t) \tan \left(\frac{p\pi}{2}
\right) \right)}$, whence, for $g\in \boL$,
  \begin{equation}
  \label{converge_dans_L}
  E_Y\left(g e^{i\frac{t}{B_n}
  S^Y_{\lfloor n m(Y) \rfloor} f_Y}\right) \to E_Y(g)E(e^{itZ})
  \end{equation}
where the random variable $Z$ is as in the statement of the
theorem (see \cite{aaronson_denker} or \cite{gouezel:stable} for
more details). We can not apply this result to $g=\phi$, since
$\phi$ is not bounded. However, $\phi$ is Lipschitz and
integrable, whence $P\phi \in \boL$ (\cite[Proposition
1.4]{aaronson_denker}). Equation \eqref{converge_dans_L} applied
to $P\phi$ gives $E_Y\left(\phi e^{i\frac{t}{B_n} S^Y_{\lfloor n
m(Y) \rfloor} f_Y\circ T_Y }\right) \to E(e^{itZ})$, since
$E_Y(P\phi)=E_Y(\phi)=1$ by Kac's Formula. Let $k(n)$ be a sequence
such that $\lfloor k(n)m(Y) \rfloor = \lfloor n m(Y) \rfloor
-1$. Since $k(n)\sim n$, the same arguments give in fact that 
$E_Y\left(\phi e^{i\frac{t}{B_n} S^Y_{\lfloor k(n)
m(Y) \rfloor} f_Y\circ T_Y }\right) \to E(e^{itZ})$, i.e.\ 
$E_Y\left(\phi e^{i\frac{t}{B_n} (S^Y_{\lfloor n
m(Y) \rfloor} f_Y  - f_Y)}\right) \to E(e^{itZ})$. The difference between
this term and $E_Y\left(\phi e^{i\frac{t}{B_n} (S^Y_{\lfloor n
m(Y) \rfloor} f_Y }\right)$ is bounded by $E_Y\left(
\phi \left|e^{-i\frac{t}{B_n}
f_Y}-1\right| \right)$, which tends to $0$ by dominated convergence. Thus,
  \begin{equation*}
  \label{eq_presque_bonne}
  E_Y\left(\phi e^{i\frac{t}{B_n}
  S^Y_{\lfloor n m(Y) \rfloor} f_Y} \right)
  \to E(e^{itZ}).
  \end{equation*}
This is \eqref{limite_mixing}. Finally, since $L$ is slowly
varying, the equation $n L(B_n)=B_n^p$ implies that 
 $\sup_{r\leq 2n} \frac{B_r}{B_n}<\infty$, $\inf_{r\geq
n} \frac{B_r}{B_n}>0$ (using for example \cite[Corollary page
274]{feller:2}).

Let $\epsilon>0$, we bound $m(M \geq \epsilon B_n)$.
  \begin{equation*}
  m(M \geq \epsilon B_n)
  \leq C (\epsilon B_n)^{-p} L( \epsilon B_n)
  =C \epsilon^{-p} B_n^{-p} L(B_n)
  \frac{L(\epsilon B_n)}{L(B_n)}.
  \end{equation*}
But $B_n^{-p}L(B_n)=\frac{1}{n}$ by definition of $B_n$, and
$\frac{L(\epsilon B_n)}{L(B_n)}$ tends to $1$ since $L$ is slowly
varying. Thus, $m(M \geq \epsilon B_n) \leq \frac{D}{n}$, which
proves \eqref{majore_fY}.

Hypothesis \ref{hypothese_3} of Theorem
\ref{thm_probabiliste_general} is satisfied for $b=1$, according
to the Birkhoff Theorem applied to $\phi-E_Y(\phi)$ (and because
$T_Y$ is ergodic, which is a consequence of the ergodicity of
$T$). Finally, the hypothesis on the distribution of $\phi$
ensures, once again by \cite{aaronson_denker}, that $\frac{S_{\lfloor
nm(Y)\rfloor}^Y
\phi -n m(Y)E_Y(\phi)}{B_n}$ converges in
distribution. Thus, \eqref{eq_hyp4'} is satisfied.
We can use Theorem
\ref{thm_probabiliste_general}, and get that $\frac{S_n f}{B_n}
\to Z$.

The proof of the second case of  Theorem \ref{thm_abstrait_markov}
is exactly the same, using \cite{aaronson_denker:central} instead
of \cite{aaronson_denker} to show the convergence in distribution
of $\frac{S_{\lfloor nm(Y)\rfloor} ^Y f_Y}{B_n}$ and $\frac{S_{\lfloor 
nm(Y) \rfloor} ^Y
\phi -n m(Y) E_Y(\phi)}{B_n}$.

In the first case ($f_Y\in L^2$), the proof is again identical
when $\phi\in L^2$, with $B_n=\sqrt{n}$, using
\cite{guivarch-hardy} (or the remarks of
\cite{aaronson_denker:central}). However, when $m(\phi>x)=x^{-p}
L(x)$, we check in a different way the hypotheses
\ref{hypothese_3} and \ref{hypothese_4} of Theorem
\ref{thm_probabiliste_general}. \cite{aaronson_denker} ensures
that, if $B'_n$ is given by
  \begin{equation}
  \label{definit_bn}
  nL(B'_n)=(B'_n)^p,
  \end{equation}
then $\frac{S_n^Y \phi -n E_Y(\phi)}{B'_n}$ converges in
distribution. Moreover, \cite[Lemma 3.4]{gouezel:stable} proves
that $Pf_Y \in \boL$, and has a vanishing integral. As $P$ has a
spectral gap on $\boL$, $P^n f_Y \to 0$ exponentially fast. In
particular, $\int f_Y\circ T_Y^n \cdot f_Y=\int (P^n f)\cdot
f=O((1-\delta)^n)$ for some $0<\delta<1$. Thus, as $f_Y\in L^2$,
\cite[Theorem 16]{vitesse_birkhoff} gives that, for every $b>1/2$,
$\frac{1}{N^b}\sum_{k=0}^{N-1}f_Y(T_Y^k) \to 0$ almost everywhere
when $N\to \infty$. In the natural extension, $\int f_Y\circ
T_Y^{-n} \cdot f_Y=\int f_Y \cdot f_Y\circ T_Y^n$ decays also
exponentially fast, whence the same argument gives that
$\frac{1}{|N|^b}\sum_{k=0}^{N-1}f_Y(T_Y^k) \to 0$ when $N\to
-\infty$. Thus, Hypothesis \ref{hypothese_3} of Theorem
\ref{thm_probabiliste_general} is satisfied for any $b>1/2$. Let
$\kappa>0$ be very small. As $L$ is slowly varying,
$L(B'_n)=O((B'_n)^\kappa)$, whence Equation \eqref{definit_bn}
gives $B'_n=O(n^{1/(p-\kappa)})$. Thus, if $b<\frac{p}{2}$, we
have $B'_n=O(B_n^{1/b})$, which implies \eqref{eq_hyp4'}.
\end{proof}

\section{Asymptotic behavior of $X_n$}
\label{section_estimee_Xn}

We return to the study of the skew product \eqref{definit_T}. To
prove limit theorems using Theorem \ref{thm_abstrait_markov}, we
will need to estimate $m(\phi_Y>n)$, which is directly related to
the speed of convergence of $X_n$ to $0$. This section will be
devoted to the proof of the following theorem:

\begin{thm}
\label{estimee_Xn_L1}
We have
  \begin{equation*}
  \left(\frac{n}{\sqrt{\ln n}}\right)^{1/\am} X_n
  \to \frac{1}{\left(2^{\am} \am^{3/2}
  \sqrt{\frac{\pi}{2\alpha''(x_0)}}\right)^{1/\am}}
  \end{equation*}
almost everywhere and in $L^1$.
\end{thm}

\begin{lem}
\label{lemme_asympt_Xn}
We have
  \begin{equation}
  E(e^{-(\alpha-\am)w})\sim \sqrt{\frac{\pi}{2\alpha''(x_0)}}
  \frac{1}{\sqrt{w}} \text{
  when }w\to\infty.
  \end{equation}
\end{lem}
\begin{proof}
Write $\beta=\alpha-\am$, and $f(b)=\Leb\{ \omega \tq
\beta(\omega) \in [0,b)\}$. In a neighborhood of $x_0$ (the unique
point where $\alpha$ takes its minimal value $\am$), $\alpha$
behaves like the parabola $\am +\frac{\alpha''(x_0)}{2}(x-x_0)^2$,
whence $f(b) \sim \sqrt{\frac{2}{\alpha''(x_0)}} \sqrt{b}$ when $b
\to 0$.

Writing $P_\beta$ for the distribution of $\beta$, an integration
by parts gives
  \begin{align*}
  E\left(e^{-(\alpha-\am)w}\right)&=\int_0^\infty e^{-b w} \dd P_\beta(b)
  =w \int_0^\infty e^{-b w} f(b) \dd b
  =\int_0^\infty e^{-u} f(u/w) \dd u
  \\ &
  =\frac{1}{\sqrt{w}} \int_0^\infty e^{-u}
  \left(\sqrt{w}f(u/w)\right) \dd u.
  \end{align*}
But $e^{-u} \left(\sqrt{w}f(u/w)\right) \to e^{-u}
\sqrt{\frac{2}{\alpha''(x_0)}}\sqrt{u}$ when $w\to \infty$. There
exists a constant $E$ such that $f(u)\leq E \sqrt{u}$ (this is
clear in a neighborhood of $0$, and elsewhere since $f$ is
bounded), whence $e^{-u} \left(\sqrt{w}f(u/w)\right) \leq E
e^{-u}\sqrt{u}$ integrable. By dominated convergence,
  \begin{equation*}
  \int_0^\infty e^{-u}
  \left(\sqrt{w}f(u/w)\right) \dd u
  \to \sqrt{\frac{2}{\alpha''(x_0)}} \int_0^\infty e^{-u} \sqrt{u}\dd u
  =\sqrt{\frac{2}{\alpha''(x_0)}} \frac{\sqrt{\pi}}{2}.
  \end{equation*}
\end{proof}

\begin{proof}[Proof of Theorem \ref{estimee_Xn_L1}]

As in Proposition \ref{estime_croissance_grossiere_Xn}, we write
  \begin{equation*}
  \frac{1}{X_n(F\omega)^{\am}}=\frac{1}{X_{n+1}(\omega)^{\am}}-\am
  2^{\am} (2X_{n+1}(\omega))^{\alpha(\omega) -\am}
  +O(X_{n+1}(\omega)^{2\alpha(\omega) -\am}).
  \end{equation*}
Proposition \ref{estime_croissance_grossiere_Xn} gives
  \begin{equation*}
  X_{n+1}(\omega)^{2\alpha(\omega)-\am} \leq
  X_{n+1}(\omega)^{\am}
 \leq
  \frac{C}{(n+1)^{\am/\aM}}
  \leq \frac{C}{\sqrt{n+1}}
  \end{equation*}
as $\am/\aM \geq 1/2$ by hypothesis. Thus,
  \begin{equation*}
  \frac{1}{X_{n+1}(\omega)^{\am}}-\frac{1}{X_n(F\omega)^{\am}}=2^{\am} \am
  (2X_{n+1}(\omega))^{\alpha-\am} +O(1/\sqrt{n}).
  \end{equation*}
Summing from $1$ to $n$, we get a constant $P$ (independent of
$\omega$) such that
  \begin{align}
  \label{minore_1/X_n}
  \frac{1}{X_n(\omega)^{\am}}\geq 2^{\am} \am
  \left[\sum_{k=1}^{n} (2X_k(F^{n-k} \omega))^{\alpha(F^{n-k}\omega)-\am}
  -P\sqrt{n} \right]
  \\
  \label{majore_1/X_n}
  \frac{1}{X_n(\omega)^{\am}}\leq 2^{\am} \am
  \left[\sum_{k=1}^n (2X_k(F^{n-k}\omega))^{\alpha(F^{n-k}\omega)-\am}
  +P\sqrt{n} \right]
  \end{align}

Equation \eqref{minore_1/X_n} and Proposition
\ref{estime_croissance_grossiere_Xn} imply that
  \begin{equation}
  \label{definit_An}
  \frac{\sqrt{\ln n}}{n} \frac{1}{2^{\am} \am X_n(\omega)^{\am}}
  \geq \frac{\sqrt{\ln n}}{n}\sum_{k=1}^n \left(\frac{2C^{-1}}{k^{1/\am}}
  \right)^{\alpha(F^{n-k} \omega)-\am}
  -P\sqrt{\frac{\ln n}{n}}=: A_n(\omega).
  \end{equation}

We first study the convergence of $A_n$. The functions $\alpha$
and $\alpha \circ F^{n-k}$ have the same distribution since $F$
preserve Lebesgue measure. Thus, by Lemma \ref{lemme_asympt_Xn},
  \begin{equation*}
  E\left(\left(\frac{2C^{-1}}{k^{1/\am}}\right)^{\alpha\circ F^{n-k}-\am}
  \right)
  \sim \sqrt{\frac{\pi}{2\alpha''(x_0)}}
  \frac{1}{\sqrt{\ln (k^{1/\am})-\ln(2C^{-1})}}
  \sim \sqrt{\frac{\pi \am}{2\alpha''(x_0)}}\frac{1}{\sqrt{\ln k}}.
  \end{equation*}
Summing, we get that
  \begin{equation}
  E(A_n)\to C_1:=\sqrt{\frac{\pi \am}{2\alpha''(x_0)}},
  \end{equation}
since $\sum_{k=2}^n \frac{1}{\sqrt{\ln k}} \sim \frac{n}{\sqrt{\ln
n}}$.

We will need $L^p$ estimates, for $p\geq 1$. To get them, we use a
result of Françoise Pène, recalled in Appendix
\ref{appendice:pene}. Let us denote by $\norm{g}$ the Lipschitz
norm of a function $g:S^1 \to \R$.

We define $f_k(\omega)=\left(\frac{2C^{-1}}{k^{1/\am}}
\right)^{\alpha(\omega)-\am}$, and $g_k=f_k-E(f_k)$. Thus,
$A_n=\frac{\sqrt{\ln n}}{n}\sum_{k=1}^n f_k \circ
F^{n-k}-P\sqrt{\frac{\ln n}{n}}$. As $g'_k=\ln
\left(\frac{2C^{-1}}{k^{1/\am}} \right)\alpha' f_k$, there exists
a constant $L$ such that, for $k\leq n$, $\norm{g_k} \leq L \ln
n$. As a consequence, Theorem \ref{thm_borne_Lp_pene} applied to
$g_k/(L \ln n)$ gives
  \begin{equation*}
  \norm{A_n-E(A_n)}_p
  = \frac{\sqrt{\ln n}}{n} L\ln n
  \norm{\sum_{k=1}^{n} g_k \circ F^{n-k} / (L \ln n)}_p
  \leq \frac{\sqrt{\ln n}}{n} L\ln n K_p \sqrt{n},
  \end{equation*}
i.e.\
  \begin{equation}
  \norm{A_n-E(A_n)}_p \leq L_p \sqrt{\frac{\ln^3 n}{n}}.
  \end{equation}

This implies in particular that $A_n$ converges almost everywhere
to $C_1$. Namely, if $\delta>0$,
  \begin{equation*}
  \Leb\{|A_n-E(A_n)|>\delta\}
  \leq \int \frac{|A_n-E(A_n)|^4}{\delta^4}
  \leq \frac{L_4}{\delta^4} \left(\frac{\ln^3 n}{n}\right)^{4/2}
  \end{equation*}
which is summable, and $E(A_n)\to C_1$.

We have
  \begin{align*}
  A_n(\omega)&
  \geq \frac{\sqrt{\ln n}}{n} \left[\sum_{k=1}^n \left( \frac{2C^{-1}}
  {k^{1/\am}}\right)^{\aM-\am} - P\sqrt{n} \right]
  \geq \frac{\sqrt{\ln n}}{n} \bigl[ K n^{2-\aM/\am} - P
  \sqrt{n}\bigr]
  \\&
  \geq K' \frac{\sqrt{\ln n}}{n} n^{2-\aM/\am}
  \end{align*}
since $\aM/\am<3/2$. Thus,
  \begin{equation}
  \norm{\frac{1}{A_n}}_\infty \leq K'' \frac{n^{\aM/\am-1}}{\sqrt{\ln
  n}}.
  \end{equation}
Note that $E(A_n)$ tends to $C_1 \not=0$, whence
$\frac{1}{E(A_n)}$ is bounded. Thus,
  \begin{align*}
  \norm{ \frac{1}{A_n} - \frac{1}{E(A_n)}}_p &
  \leq \norm{\frac{1}{A_n}}_\infty \frac{1}{E(A_n)}
  \norm{A_n-E(A_n)}_p
  \leq K''' \frac{n^{\aM/\am-1}}{\sqrt{\ln
  n}} L_p \sqrt{\frac{\ln^3 n}{n}}
  \\&
  = M_p \frac{\ln n}{n^\kappa}
  \end{align*}
where $\kappa=\frac{3}{2}-\frac{\aM}{\am}>0$. In particular,
$\frac{1}{A_n}$ tends to $\frac{1}{C_1}$ in every $L^p$. Equation
\eqref{definit_An} shows that
  \begin{equation}
  \label{un_petit_label}
  \left(\frac{n}{\sqrt{\ln n}}\right)^{1/\am} X_n \leq \frac{1}
  {(2^{\am}\am A_n )^{1/\am}}.
  \end{equation}
The right hand side tends to
  \begin{equation}
  C_2:=\frac{1}{\left(2^{\am}\am^{3/2}
  \sqrt{\frac{\pi}{2\alpha''(x_0)}}\right)^{1/\am}}
  \end{equation}
in every $L^p$, and in particular in $L^1$. Thus,
  \begin{equation}
  \label{eq_esperance}
  \varlimsup E\left(\left( \frac{n}{\sqrt{\ln n}}\right)^{1/\am} X_n
  \right) \leq C_2.
  \end{equation}
Moreover, $A_n$ converges almost everywhere to $C_1$, whence
\eqref{un_petit_label} yields that, almost everywhere,
  \begin{equation}
  \label{eq_limsup}
  \varlimsup \left( \frac{n}{\sqrt{\ln n}}\right)^{1/\am}
  X_n(\omega)
  \leq C_2.
  \end{equation}

Set $Q=\sup_n \left(\frac{1}{E(A_n)}\right)+1$, we estimate
$\Leb\left\{ \frac{1}{A_n} \geq Q\right\}$. If $p\geq 1$,
  \begin{equation*}
  \Leb\left\{ \frac{1}{A_n} \geq Q\right\}
  \leq \Leb \left\{ \left|\frac{1}{A_n} -\frac{1}{E(A_n)} \right|
  \geq 1 \right\}
  \leq E\left( \left|\frac{1}{A_n} -\frac{1}{E(A_n)} \right| ^p
  \right)
  \leq \left( M_p \frac{\ln n}{n^\kappa} \right)^p.
  \end{equation*}
In particular, choosing $p$ large enough gives
  \begin{equation*}
  \Leb\left\{ \frac{1}{A_n} \geq Q\right\}\leq  \frac{M}{n^5}.
  \end{equation*}
Setting $Q'=\frac{Q}{2^{\am} \am}$, \eqref{un_petit_label} thus
yields that
  \begin{equation}
  \Leb\left\{ X_n \geq \left(\frac{Q'\sqrt{\ln n}}{n}\right)^{1/\am}
  \right\} \leq \frac{M}{n^5}.
  \end{equation}
Consequently, $U_n:=\left\{ \omega \tq \exists \sqrt{n} \leq k\leq
n \text{ with } X_k(F^{n-k} \omega) \geq \left(\frac{Q'\sqrt{\ln
k}}{k}\right)^{1/\am} \right\}$ has a measure at most
$\sum_{\sqrt{n}}^n \frac{M}{k^5} \leq \frac{M'}{n^2}$ (since
$\Leb$ is invariant under $F^{n-k}$). Finally, Borel-Cantelli
ensures that there is a full measure subset of $S^1$ on which
$\omega \not \in U_n$ for large enough $n$.

Set
  \begin{equation*}
  A'_n(\omega)=\frac{\sqrt{\ln n}}{n} \left[\sum_{k=1}^n
  \left(\frac{2(Q'\sqrt{\ln
  k})^{1/\am}}{k^{1/\am}}\right)^{\alpha(F^{n-k} \omega) -\am}
  + (P+1)\sqrt{n}\right].
  \end{equation*}
As for $A_n$, we show that $A'_n \to C_1$ in every $L^p$ and
almost everywhere.

Let $\omega$ be such that $\omega \not \in U_n$ for large enough
$n$, and $A'_n(\omega) \to C_1$ (these properties are true almost
everywhere). Then, for large enough  $n$, Equation
\eqref{majore_1/X_n} and the fact that $X_k(F^{n-k} \omega) \leq
\left(\frac{Q'\sqrt{\ln k}}{k}\right)^{1/\am}$ for $\sqrt{n}\leq
k\leq n$, yield that
  \begin{align*}
  \frac{1}{2^{\am} \am X_n(\omega)^{\am}}&
  \leq  \left[\sum_{k=1}^{\sqrt{n}} 1+
  \sum_{k=\sqrt{n}}^n \left(\frac{2(Q'\sqrt{\ln
  k})^{1/\am}}{k^{1/\am}}\right)^{\alpha(F^{n-k}\omega) -\am}
  +P \sqrt{n} \right]
  \\&
  \leq \frac{n}{\sqrt{\ln n}} A'_n(\omega)\sim \frac{n}{\sqrt{\ln n}} C_1.
  \end{align*}
Thus,
  \begin{equation}
  \label{eq_liminf}
  \varliminf \left( \frac{n}{\sqrt{\ln n}}\right)^{1/\am}
  X_n(\omega)
  \geq C_2.
  \end{equation}

Equations \eqref{eq_limsup} and \eqref{eq_liminf} prove that
$\left( \frac{n}{\sqrt{\ln n}}\right)^{1/\am} X_n$ tends almost
everywhere to $C_2$. We get the convergence in $L^1$ from the
inequality \eqref{eq_esperance} and the following elementary
lemma.
\end{proof}

\begin{lem}
Let $f_n$ be nonnegative functions on a probability space, with
$f_n \to f$ almost everywhere, and $\varlimsup E(f_n) \leq
E(f)<\infty$. Then $f_n \to f$ in $L^1$.
\end{lem}

\section{Limit theorems}

\label{section:limite}

Set
  \begin{equation}
  \label{definit_A}
  A=\frac{1}{4\left( \am^{3/2}
  \sqrt{\frac{\pi}{2\alpha''(x_0)}}\right)^{1/\am}}\int_{S^1\times \{1/2\}}
  h\dLeb,
  \end{equation}
where $h$ is the density of $m$ with respect to $\Leb$.

In this section, we prove the following theorem:

\begin{thm}
\label{enonce_theoreme_limite}
Let $f$ be a Hölder function on $S^1\times [0,1]$, with $\int f\dd
m=0$. Write $c=\int_{S^1\times\{0\}} f \dLeb$. Then
  \begin{itemize}
  \item If $\am<1/2$, there exists $\sigma^2 \geq 0$ such that
$\frac{1}{\sqrt{n}} S_n f \to \boN(0,\sigma^2)$.
  \item If $\am=1/2$ and $c\not=0$, then
$\frac{S_n f}{\sqrt{ \frac{c^2 A }{4} n (\ln n)^2}} \to
\boN(0,1)$.
  \item If $1/2<\am<1$ and $c\not =0$, then
$\frac{S_n f}{n^{\am}\sqrt{\am \ln n}} \to Z$, where the random
variable $Z$ has a characteristic function given by
  \begin{equation*}
  E(e^{itZ})
  =e^{- A
  |c|^{1/\am} \Gamma(1-1/\am)\cos\left(\frac{\pi}{2\am}\right) |t|^{1/\am}
  \left( 1-i\sgn(ct) \tan
  \left(\frac{\pi}{2\am} \right) \right)}
  \end{equation*}
  \item If $1/2\leq \am<1$ and $c=0$, assume also that there exists $\gamma
>0$ such that $|f(\omega,x)-f(\omega,0)|\leq Cx^\gamma$, with
$\gamma>\aM\left(1-\frac{1}{2\am}\right)$. Then there exists
$\sigma^2 \geq 0$ such that $\frac{1}{\sqrt{n}} S_n f \to
\boN(0,\sigma^2)$.
  \end{itemize}
\end{thm}

The random variable $Z$ in the third case has a stable
distribution of exponent $1/\am$ and parameters $A|c|^{1/\am}
\Gamma(1-1/\am)\cos\left(\frac{\pi}{2\am}\right)$ and $\sgn(c)$.

To prove this theorem, we will use Theorem
\ref{thm_abstrait_markov}. For this, we need a control of
$m(\phi_Y>n)$ which comes from the asymptotic behavior of $X_n$
proved in Theorem \ref{estimee_Xn_L1}. It will also be necessary
to estimate $m(f_Y>x)$, through the study of the integrability of
$f_Y$ (Lemmas \ref{L2_am_petit} and \ref{lemme_dans_Lp}).

In the rest of this section, $f$ will be a Hölder function on
$S^1\times [0,1]$, fixed once and for all. Recall that
$f_Y(y)=\sum_{k=0}^{\phi_Y(y)-1} f(T^k y)$, where $\phi_Y$ is the
first return time to $Y=S^1\times(1/2,1]$.

\subsection{Estimates on measures}
\begin{lem}
\label{controle_mesure_phi}
We have
  \begin{equation*}
  m(\phi_Y>n) \sim \left(\frac{\sqrt{\ln n}}{n}\right)^{1/\am} A
  \end{equation*}
where $A$ is given by \eqref{definit_A}.
\end{lem}
\begin{proof}
We have
  \begin{align*}
  m(\phi_Y> n)&
  =\int_{S^1} \int_{1/2}^{Y_{n+1}(\omega)} h(\omega,u)\dd u \dd
  \omega
  =\int_{S^1} \int_0^{X_{n}(F\omega)/2} h(\omega,1/2+u)\dd u \dd\omega
  \\&
  =\int_{S^1} \frac{X_{n}(F\omega)}{2} h(\omega,1/2)\dd\omega
  +\int_{S^1} \int_0^{X_{n}(F\omega)/2} \bigl[h(\omega,1/2+u)-h(\omega,1/2)
  \bigr] \dd u\dd \omega
  \\&
  =I+II.
  \end{align*}
As $\left(\frac{n}{\sqrt{\ln n}}\right)^{1/\am} X_n(F\omega) \to
\frac{1}{\left(2^{\am} \am^{3/2}
  \sqrt{\frac{\pi}{2\alpha''(x_0)}}\right)^{1/\am}}$ in
$L^1$ and almost everywhere (Theorem \ref{estimee_Xn_L1}) and
$h(\omega,1/2)$ is bounded, we get that $I\sim
\left(\frac{\sqrt{\ln n}}{n}\right)^{1/\am} A$. Moreover, for
large enough $n$, $|h(\omega,1/2+u)-h(\omega,1/2)|\leq \epsilon$,
whence $II=o\left(\frac{\sqrt{\ln n}}{n}\right)^{1/\am}$.
\end{proof}

\begin{lem}
\label{L2_am_petit}
If $\am<1/2$, then $f_Y \in L^2(Y, \de m)$.
\end{lem}
\begin{proof}
We have
  \begin{align*}
  \int f_Y^2 \dd m &
  \leq C\sum_n m(\phi_Y=n)  n^2
  =C \sum \bigl( m(\phi_Y>n-1)-m(\phi_Y>n)\bigr)n^2
  \\&
  \leq C\sum m(\phi_Y>n) n
  \end{align*}
which is summable since $m(\phi_Y>n)\sim A \left(\frac{\sqrt{\ln
n}}{n}\right)^{1/\am}$ with $1/\am>2$.
\end{proof}

\begin{lem}
\label{lemme_dans_Lp}
Assume that $\int_{S^1\times \{0\}} f=0$. Let $\aM>\gamma>0$ be
such that $|f(\omega,x)-f(\omega,0)| \leq C x^\gamma$. If $1<p<
\min\left( \frac{2}{\am}, \frac{1}{\am(1-\gamma/\aM)} \right)$,
then $f_Y \in L^p(Y,\de m)$.
\end{lem}
\begin{proof}
As $h$ is bounded on $Y$, it is sufficient to prove that $f_Y \in
L^p(Y,\dLeb)$.

Assume first that $f \equiv 0$ on $S^1 \times \{0\}$. Then, if
$\xx=(\omega,x)$ satisfies $\phi_Y(\xx)=n$, we have
$f_Y(\xx)=\sum_0^{n-1} f(T^k \xx)$. If $k \geq 1$,
$T_\omega^k(x)\leq X_{n-k}(F^k \omega) \leq
\frac{C}{(n-k)^{1/\aM}}$, whence $|f(T^k \xx)|\leq
\frac{C}{(n-k)^{\gamma/\aM}}$, and a summation yields that
$|f_Y(\xx)|\leq C n^{1-\gamma/\aM}$. Thus,
  \begin{align*}
  \int |f_Y|^p &
  \leq C \sum m(\phi_Y=n) n^{p(1-\gamma/\aM)}
  \\&
  \leq C \sum m(\phi_Y>n) n^{p(1-\gamma/\aM)-1}.
  \end{align*}
As $m(\phi_Y>n) \sim A \left(\frac{\sqrt{\ln n}}{n}
\right)^{1/\am}$, this last series is summable as soon as
  \begin{equation*}
  -\frac{1}{\am}+p\left(1-\frac{\gamma}{\aM}\right)-1 <-1,
  \end{equation*}
which is the case by assumption on $p$.

Assume now that $f$ has a vanishing integral on $S^1$. Let
$g(\omega,x)=f(\omega,0)$. The function $f-g$ vanishes on $S^1
\times \{0\}$, whence $f_Y-g_Y \in L^p$ according to the first
part of this proof. Consequently, it is sufficient to prove that
$g_Y \in L^p$. Write $\chi(\omega)=f(\omega,0)$ and $S_n
\chi(\omega)=\sum_{k=0}^{n-1} \chi(F^k \omega)$: then
$g_Y(\omega,x)=S_{\phi_Y(\omega,x)} \chi(\omega)$.

Let $M_n \chi(\omega)=\max_{k\leq n}|S_k \chi(\omega)|$. Let
$\delta>0$, and $l=\frac{1+\delta}{\delta}$, so that
$\frac{1}{l}+\frac{1}{1+\delta}=1$. We have
  \begin{align*}
  \int |g_Y|^p &
  =\sum_{n=0}^{\infty} \int_{S^1}
  \int_{1/2+X_{n}(F\omega)/2}^{1/2+X_{n-1}(F\omega)/2} \bigl|S_n
  \chi(\omega)\bigr|^p \dd u \dd \omega
  \\&
  \leq \sum_{k=1}^{\infty} \int_{S^1}
  \int_{1/2+X_{2^k}(F\omega)/2}^{1/2+X_{2^{k-1}}(F\omega)/2} |M_{2^k}
  \chi(\omega)|^p \dd u \dd \omega
  \\&
  \leq  \sum_{k=1}^\infty \int_{S^1} X_{2^{k-1}}(F\omega)
  |M_{2^k}\chi(\omega)|^p \dd
  \omega
  \leq  \sum_{k=1}^\infty \norm{X_{2^{k-1}}\circ F}_{1+\delta}
  \norm{M_{2^k}\chi}_{lp}^p,
  \end{align*}
where the last inequality is Hölder inequality. If $\delta$ is
small enough, $lp>2$, whence Corollary \ref{ineg_max} yields that
$\norm{M_{2^k}\chi}_{lp} \leq C k^{\frac{lp-1}{lp}}\sqrt{2^k}$.
Moreover,
  \begin{equation*}
  \norm{X_{2^{k-1}}\circ F}_{1+\delta}
  =\norm{X_{2^{k-1}}}_{1+\delta}
  \leq \left( \int X_{2^{k-1}} \right)^{1/(1+\delta)}
  \sim C \left( \frac{\sqrt{\ln(2^{k-1})}}{2^{k-1}}
  \right)^{\frac{1}{(1+\delta) \am}}
  \end{equation*}
by Theorem \ref{estimee_Xn_L1}. Thus, $\int|g_Y|^p<\infty$ if
$\frac{1}{(1+\delta)\am} > \frac{p}{2}$, and it is possible to
choose $\delta$ such that this inequality is true, since
$\frac{1}{\am}>\frac{p}{2}$ by hypothesis.
\end{proof}

\subsection{Proof of Theorem \ref{enonce_theoreme_limite}}

To apply Theorem \ref{thm_abstrait_markov}, we first check the
condition \eqref{condition_markov_abstrait}. Let $\theta$ be the
Hölder exponent of $f$. We will work with the distance $d_{
\lambda^{-\theta}}= \lambda^{-\theta s(x,y)}$. For this distance,
$T_Y$ is a Gibbs-Markov map.

\emph{Fact: if $f$ is $\theta$-Hölder on $S^1\times [0,1]$, then
  \begin{equation*}
  \sum m[A_{s,n}] D f_Y(A_{s,n}) <\infty.
  \end{equation*}
}

Recall that $D f_Y(A_{s,n})$ (defined in Theorem
\ref{thm_abstrait_markov}) is the best Lipschitz constant of $f_Y$
on $A_{s,n}$, here for the distance $d_{ \lambda^{-\theta}}$.

\begin{proof}
Take $(\omega_1,x_1)$ and $(\omega_2,x_2)\in A_{s,n}$ with for
example $x_2\geq x_1$. This implies that $x_1\in J_n^+(\omega_2)$
and that, for $0\leq k\leq n$,
$d(T^k(\omega_1,x_1),T^k(\omega_2,x_2)) \leq (1+D) |F^k \omega_1
-F^k \omega_2|$ (see the beginning of the proof of Proposition
\ref{prop_distortion_bornee}). Moreover,
$d(T^k(\omega_2,x_1),T^k(\omega_2,x_2)) \leq
d(T^n(\omega_2,x_1),T^n(\omega_2,x_2))$ (since, if $\omega$ is
fixed, the map $T_{\alpha(\omega)}$ is expanding).

Thus, for $0\leq k\leq n$,
  \begin{align*}
  d(T^k(\omega_1,x_1),T^k(\omega_2,x_2))&
  \leq
  d(T^k(\omega_1,x_1),T^k(\omega_2,x_1))
  +d(T^k(\omega_2,x_1),T^k(\omega_2,x_2))
  \\&
  \leq
  (1+D)|F^k \omega_1 -F^k \omega_2|
  +d(T^n(\omega_2,x_1),T^n(\omega_2,x_2))
  \\&
  \leq
  (1+D)|F^n \omega_1 -F^n \omega_2|
  +d(T^n(\omega_1,x_1),T^n(\omega_2,x_1))
  \\& \hphantom{=\ }
  +d(T^n(\omega_1,x_1),T^n(\omega_2,x_2))
  \\&
  \leq
  (1+D)|F^n \omega_1 -F^n \omega_2|
  +(1+D)|F^n \omega_1 -F^n \omega_2|
  \\& \hphantom{=\ }
  +d(T^n(\omega_1,x_1),T^n(\omega_2,x_2))
  \\&
  \leq
  (3+2D)d(T^n(\omega_1,x_1),T^n(\omega_2,x_2)).
  \end{align*}
We deduce that
  \begin{align*}
  |f_Y(\omega_1,x_1)-f_Y(\omega_2,x_2)|&
  \leq \sum_{k=0}^{n-1}
  |f(T^k(\omega_1,x_1))-f(T^k(\omega_2,x_2))|
  \\&
  \leq \sum_{k=0}^{n-1}C
  d(T^k(\omega_1,x_1),T^k(\omega_2,x_2))^\theta
  \\&
  \leq C' n d(T^n(\omega_1,x_1),T^n(\omega_2,x_2))^\theta.
  \end{align*}
As $T_Y$ is expanding for the distance $d'$ (defined in
\eqref{definit_d'}, and equivalent to $d$), we get
  \begin{equation*}
  d(T^n(\omega_1,x_1),T^n(\omega_2,x_2)) \leq C
  d_{\lambda^{-1}}(T^n(\omega_1,x_1),T^n(\omega_2,x_2))=C \lambda
  d_{\lambda^{-1}}((\omega_1,x_1),(\omega_2,x_2)),
  \end{equation*}
whence $d(T^n(\omega_1,x_1),T^n(\omega_2,x_2))^\theta \leq C
d_{\lambda^{-\theta}}((\omega_1,x_1),(\omega_2,x_2))$.

Thus, $Df_Y(A_{s,n}) \leq C n$, and
  \begin{equation*}
  \sum m(A_{s,n}) D f_Y(A_{s,n})
  \leq C \sum m(\phi_Y=n) n
  = C < +\infty,
  \end{equation*}
by Kac's Formula. \qedfact
\end{proof}

\begin{proof}[Proof of Theorem \ref{enonce_theoreme_limite}]

In the case $\am<1/2$, Lemma \ref{L2_am_petit} gives that $f_Y \in
L^2$. Moreover, $|f|_Y\in L^2$ for the same reason, and $\phi \in
L^2$ (since $\phi=g_Y$ for $g\equiv 1$, whence Lemma
\ref{L2_am_petit} applies also). We have already checked the
condition \eqref{condition_markov_abstrait}, so we can apply (the
first case of) Theorem \ref{thm_abstrait_markov}. This yields the
central limit theorem for $f$.

The second and third cases are analogous. Let us prove for example
the third one, i.e.\ $1/2<\am<1$ and $c\not=0$. Assume for example
$c>0$. We estimate $m(f_Y>x)$.

\emph{Fact: $m(f_Y>x) \sim \left(\frac{c \sqrt{\ln
x}}{x}\right)^{1/\am} A$ and $m(f_Y<-x)=o\left(\frac{\sqrt{\ln
x}}{x}\right)^{1/\am}$.}
\begin{proof}
We prove the estimate on $m(f_Y>x)$, the other one being similar.

Let $g\equiv c$ on $S^1 \times [0,1]$. Then $g_Y=nc$ on
$[\phi_Y=n]$, which implies that $m(g_Y>nc)=m(\phi_Y>n)\sim
\left(\frac{\sqrt{\ln n}}{n}\right)^{1/\am} A$ by Lemma
\ref{controle_mesure_phi}.

In the general case, consider $j=f-g$, and let us prove that
$m(|j_Y|>x)=o\left(\frac{\sqrt{\ln x}}{x}\right)^{1/\am}$. As
$f_Y=g_Y+j_Y$, it will give
  \begin{equation*}
  m(g_Y>x(1+\epsilon)) -m(|j_Y|>x\epsilon) \leq m(f_Y>x) \leq
  m(g_Y>x(1-\epsilon))+m(|j_Y|>x\epsilon),
  \end{equation*}
which gives the conclusion.

Let $\gamma>0$ with $\gamma<\min(\theta,\aM)$ (where $\theta$ is
the Hölder coefficient of $f$). Lemma \ref{lemme_dans_Lp} gives
that $j_Y\in L^p$ if
$p<\min\left(\frac{2}{\am},\frac{1}{\am(1-\gamma/\aM)} \right)$.
We can in particular choose $p>1/\am$. Then $m(|j_Y|>x) \leq \int
\left(\frac{|j_Y|}{x}\right)^p =O(x^{-p})$, which concludes the
proof of the fact. \qedfact
\end{proof}

The same fact holds for $\phi_Y$ and $|f|_Y$, with the same proof,
whence we are in the third case of Theorem
\ref{thm_abstrait_markov}. This gives the desired result.

Assume finally that $\frac{1}{2}\leq \am<1$ and that $c=0$. Under
the hypotheses of the theorem, we can apply Lemma
\ref{lemme_dans_Lp} with $p=2$, and get that $f_Y \in L^2$. The
proof of this lemma shows in fact that the function $M$ (defined
in Theorem \ref{thm_abstrait_markov}) is also in $L^2$. Finally,
Lemma \ref{controle_mesure_phi} shows that $m[\phi_Y>x] \sim
\left(\frac{\sqrt{\ln x}}{x}\right)^{1/\am}A$. We have checked all
the hypotheses of the first case of Theorem
\ref{thm_abstrait_markov}.
\end{proof}

\appendix

\section{Induced maps and limit theorems}

\label{appendice:loi_stable}

The aim of this section is to prove very general results stating
that, if a function satisfies a limit theorem for an induced map,
it also satisfies one for the initial map. Similar theorems have
been proved in \cite{gouezel:stable}, by spectral methods. We will
describe here a more elementary method, essentially due to
Melbourne and Török for flows (\cite{melbourne_torok}).

If $Y$ is a subset of a probability space $(X,m)$, $T:X\to X$, and
$T_Y$ is the induced map on $Y$, we will write $S_n^Y
g=\sum_{k=0}^{n-1} g\circ T_Y^k$: this is the Birkhoff sum of $g$,
for the transformation $T_Y$. We will also write
$E_Y(g)=\frac{\int_Y g}{m[Y]}$. Finally, for $t\in \R$, $\lfloor t
\rfloor$ denotes the integer part of $t$.

\begin{thm}
\label{thm_probabiliste_general}
Let $T:X\to X$ be an ergodic endomorphism of a probability space
$(X,m)$, and $f:X\to \R$ an integrable function with vanishing
integral. Let $Y\subset X$ have positive measure. For $y\in Y$,
write $\phi(y)=\inf\{n>0 \tq T^n(y)\in Y\}$ and
$f_Y(y)=\sum_{k=0}^{\phi(y)-1} f(T^k y)$, and $M(y)=\max_{1\leq
k\leq \phi(y)} \left|\sum_{j=0}^{k-1} f(T^j y) \right|$.

We assume the following properties:
\begin{enumerate}
\item
There exists a sequence $B_n\to +\infty$, with $\sup_{r \leq 2n}
\frac{B_r}{B_n} < \infty$ and $\inf_{r \geq n} \frac{B_r}{B_n}>0$,
such that $(f_Y, \phi)$ satisfies a mixing limit theorem for the
normalization $B_n$: there exists a random variable $Z$ such that,
for every $t\in \R$,
\begin{equation}
  \label{limite_mixing}
  E_Y\left(\phi e^{it \frac{S^Y_{\lfloor n m(Y) \rfloor}  f_Y}{B_n}}\right)
  \to E_Y(\phi) E\left(e^{itZ}\right).
  \end{equation}
\item
For every $\epsilon>0$, there exists $C$ such that, for any $n\in
\N^*$,
\begin{equation}
  \label{majore_fY}
  m\{y\in Y \tq M(y) \geq \epsilon B_n \} \leq \frac{C}{n}.
  \end{equation}
\item
\label{hypothese_3}
There exists $b>0$ such that, in the natural extension of $T_Y$,
$\frac{1}{N^b} \sum_0^{N-1} f_Y(T_Y^k)$ tends almost everywhere to
$0$ when $N \to \pm \infty$.
\item
\label{hypothese_4}
For every $\epsilon>0$, there exists $A>0$ and $N_0$ such that,
for every $n\geq N_0$,
  \begin{equation}
  \label{eq_hyp4}
  m\left\{ y\in Y \tq \left| \frac{S_n^Y \phi - n
  E_Y(\phi)}{B_n^{1/b}} \right| \geq A \right\} \leq \epsilon.
  \end{equation}
\end{enumerate}

Then the function $f$ satisfies also a limit theorem:
\begin{equation*}
  E\left(e^{it \frac{S_n f}{B_n}}\right)
  \to E(e^{it Z}),
  \end{equation*}
i.e.\ $\frac{S_n f}{B_n}$ tends in distribution to $Z$.
\end{thm}

The hypotheses of the theorem are tailor-made so that the
following proof works, but they are in fact often satisfied in
natural cases. Let us comment on these 4 hypotheses:
\begin{enumerate}
\item
The convergence \eqref{limite_mixing} is very often satisfied when
$f_Y$ satisfies a limit theorem. Namely, the martingale proofs or
spectral proofs of limit theorems automatically give this kind of
convergence.
\item
If $Z_0,Z_1,\ldots$ are independent identically random variables
such that $\frac{\sum_{0}^{n-1} Z_k}{B_n}$ converges in
distribution to a nontrivial limit, then for all $\epsilon>0$,
there exists $C$ such that $P(|Z_0| \geq \epsilon B_n) \leq
\frac{C}{n}$: this is a consequence of the classification of the
stable laws, see \cite{feller:2}.

Here, we are not in the independent setting, and there is no such
classification. However, the same kind of results holds very
often: usually, it is not hard to check in practical cases that
$m(|f_Y(x)| \geq \epsilon B_n) \leq \frac{C}{n}$, since $f_Y$
satisfies a limit theorem by the first assumption. Set
$|f|_Y(y)=\sum_{j=0}^{\phi(y)-1} |f(T^j y)|$. As $|f|_Y$ and $f_Y$
have more or less the same distribution, $|f|_Y$ satisfies also
often
  \begin{equation}
  \label{majore_fY'}
  \tag{\ref{majore_fY}'}
  m\{ y\in Y \tq |f|_Y(y) \geq \epsilon B_n \} \leq \frac{C}{n}.
  \end{equation}
Since $M \leq f_Y$, \eqref{majore_fY'} implies \eqref{majore_fY}.
Thus, it will often be sufficient to check \eqref{majore_fY'}.
However, \eqref{majore_fY} is sometimes strictly weaker than
\eqref{majore_fY'}, because of cancellations, which is why we have
stated the theorem with \eqref{majore_fY}.
\item
The natural extension is useful so that we can let $N$ tend to
$-\infty$, and consider $T_Y^{-1}$ in the proof. Generally,
Birkhoff's Theorem yields that this assumption is satisfied for
$b=1$. This is often sufficient. However, sometimes, it is
important to have better estimates. It is then possible to use
\cite[Theorem 16]{vitesse_birkhoff}, for example: this theorem
ensures that, if the correlations of $f_Y\in L^2$ decay at least
as $O(1/n)$, then the hypothesis is satisfied for any $b>1/2$ (for
$N \to -\infty$, use the fact that $\int f_Y \cdot f_Y \circ T_Y^n
= \int f_Y\circ T_Y^{-n} \cdot f_Y$, and apply the result to
$T_Y^{-1}$).
\item
The fourth assumption is weaker than
  \begin{equation}
  \label{eq_hyp4'}
  \tag{\ref{eq_hyp4}'}
  \exists B'_n=O(B_n^{1/b}) \text{ such that }
  \frac{S_n^Y \phi - nE_Y(\phi)}{B'_n}\text{ converges in
  distribution.}
  \end{equation}
Moreover, $\phi$ is often simpler than $f_Y$. Since $f_Y$
satisfies a limit theorem (this is more or less the first
hypothesis), this is also often the case of $\phi$, which implies
\eqref{eq_hyp4'}. Thus,
 \eqref{eq_hyp4'} -- and hence  \eqref{eq_hyp4} --
are satisfied quite generally.
\end{enumerate}

\begin{proof}[Proof of Theorem \ref{thm_probabiliste_general}]
Without loss of generality, we can work in a tower, i.e.\ assume
that $X=\{ (y,i) \tq y\in Y, i\in\{0,\ldots,\phi(y)-1\} \}$ and
that, for $i<\phi(y)-1$, $T(y,i)=(y,i+1)$, while
$T(y,\phi(y)-1)=(T_Y(y),0)$. Namely, it is possible to build an
extension of $X$ satisfying these properties, and it is equivalent
to prove a limit theorem in $X$ or in this extension (see for
example \cite[Section 4.1]{gouezel:stable}). Note that
$E_Y(\phi)=1/m(Y)$ by Kac's Formula. Let $\pi$ be the projection
from $X$ to $Y$, given by $\pi(y,i)=y$.

In this proof, we will write $S_t f(x)$, even when $t$ is not an
integer, for $S_{\lfloor t \rfloor }f(x)$. In the same way, $T^t$
should be understood as $T^{\lfloor t \rfloor}$. We also extend
$B_n$ to $\R_+$, setting $B_t:=B_{\lfloor t \rfloor}$.

As $T$ is ergodic, $T_Y$ is also ergodic (\cite[Proposition
1.5.2]{aaronson:book}). Birkhoff's Theorem gives that
  \begin{equation}
  \label{asymp_phin}
  S_n^Y \phi = \frac{n}{m(Y)}  +o(n)
  \end{equation}
almost everywhere on $Y$. For $y\in Y$ and $N\in \N$, let $n(y,N)$
be the greatest integer $n$ such that $S_n^Y \phi(y) < N$. If $y$
is such that $S_n^Y \phi(y)=\frac{n}{m(Y)}+o(n)$ (which is true
almost everywhere), then $n(y,N)$ is
finite for every $N$, and $\frac{n(y,N)}{m(Y)} \sim N$, i.e.\
\begin{equation}
  \label{renouvellement_basique}
  \frac{n(y,N)}{N m(Y)} \to 1.
  \end{equation}

Since $\int_X e^{it (S_N^Y f_Y)\circ \pi }=\int_Y \phi e^{it S_N^Y
f_Y}$, \eqref{limite_mixing} yields that
  \begin{equation}
  \label{Thm_limite_SY}
  \frac{(S^Y_{Nm(Y)} f_Y)\circ \pi}{B_N} \to Z
  \end{equation}
in distribution on $X$. The idea of the proof will be to see that
$(S^Y_{Nm(Y)}f_Y)\circ \pi$ and $S_N f$ are close (this is not
surprising, since one iteration of $T_Y$ corresponds roughly to
$1/m(Y)$ iterations of $T$). This will give that $\frac{S_N
f}{B_N}$ tends to $Z$.

We write
  \begin{align*}
  S_N f(y,i)=
  \left(S_N f(y,i)-S_N f(y,0) \right)
  &+ \left(S_N f(y,0)-S^Y_{n(y,N)}f_Y(y) \right)
  \\&
  +\left( S^Y_{n(y,N)}f_Y(y) - S^Y_{Nm(Y)}f_Y(y) \right)
  +S^Y_{ Nm(Y)}f_Y(y).
  \end{align*}
The last term, equal to $\bigl(S^Y_{N m(Y)} f_Y\bigr)\circ \pi$,
satisfies a limit theorem by \eqref{Thm_limite_SY}. To conclude
the proof, we will see that the three other terms, divided by $B_N$, 
tend to $0$ in probability.

The second and third terms depend only on $y$. Thus, the following
lemma will be useful to prove that they tend to $0$ on $X$:
\begin{lem}
\label{proba_sur_Y_donne_proba_sur_X}
Let $f_n$ be a sequence of functions on $Y$, tending to $0$ in
probability on $Y$. Then $f_n \circ \pi$ tends to $0$ in
probability on $X$.
\end{lem}
\begin{proof}
Take $\epsilon>0$. As $f_n \to 0$ in probability, the measure of
$E_n:=\{ y\in Y \tq |f_n(y)|\geq \epsilon\}$ tends to $0$. As
$\phi \in L^1$, dominated convergence yields that $\int_{E_n} \phi
\to 0$, i.e. the measure of $\pi^{-1}(E_n)$ tends to $0$. But
$\pi^{-1}(E_n)$ is exactly the set where $|f_n \circ \pi|\geq
\epsilon$.
\end{proof}

\emph{Fact: $\frac{1}{B_N} \left(S_N f(y,i)-S_N f(y,0) \right)$
tends to $0$ in probability on $X$.}
\begin{proof}
Set $V_N(y)=\sum_{i=0}^{\phi(y)-1} |f \circ T^N(y,i)|$ on $Y$.
Then $\norm{V_N}_{L^1(Y)} =\norm{ f\circ T^N}_{L^1(X)}
=\norm{f}_{L^1(X)}$ since $T$ preserves the measure. Thus,
$V_N/B_{N}$ tends to $0$ in $L^1(Y)$, and in probability. Lemma
\ref{proba_sur_Y_donne_proba_sur_X} yields that $\frac{1}{B_N}
V_N\circ \pi$ tends to $0$ in probability on $X$.

As $S_N f(y,i)-S_N f(y,0) =\sum_N^{N+i-1} f(T^k(y,0))
-\sum_0^{i-1} f(T^k(y,0))$, we get $|S_N f(y,i)-S_N f(y,0)| \leq
V_N(y)+V_0(y)$. Thus, $\frac{1}{B_N} \left(S_N f(y,i)-S_N f(y,0)
\right)$ is bounded by a function going to $0$ in probability.
\qedfact
\end{proof}

\emph{Fact:  $\frac{1}{B_N} \left(S_N f(y,0)-S_{n(y,N)}^Y
f_Y(y)\right)$ tends to $0$ in probability on $X$.}
\begin{proof}
By Lemma  \ref{proba_sur_Y_donne_proba_sur_X}, it is sufficient to
prove it on $Y$. We have
  \begin{equation*}
  \left|S_N f(y,0)-S_{n(y,N)}^Y f_Y(y)\right|
  =
  \left| \sum_{S_{n(y,N)}^Y \phi(y)}^{N-1} f\circ T^k(y,0) \right|
  \leq M\left(T_Y^{n(y,N)}y\right).
  \end{equation*}
Let $a>0$ be very small, we show that $m\left\{ y \tq
M\left(T_Y^{n(y,N)}y\right) \geq a B_N \right\} \to 0$.

Let $\epsilon>0$. Let $C$ be such that $m( M(y) \geq a B_n) \leq
\frac{C}{n}$, by \eqref{majore_fY}. Set
$\delta=\frac{\epsilon}{2Cm(Y)}$. By
\eqref{renouvellement_basique}, for large enough $N$,
  \begin{equation*}
  m\left\{\left|\frac{n(y,N)}{m(Y) N} - 1 \right|
  \geq \delta\right\}
  \leq \epsilon.
  \end{equation*}
When $\left|\frac{n(y,N)}{m(Y) N} - 1 \right| \leq \delta$, the
fact that $M\left(T_Y^{n(y,N)}y\right) \geq aB_N$ implies that
there exists $n \in [(1-\delta) m(Y)N, (1+\delta)m(Y)N]$ such that
$M(T_Y^n y) \geq aB_N$. Thus,
  \begin{equation*}
  m\left\{ y \tq M\left(T_Y^{n(y,N)}y\right) \geq a B_N \right\}
  \leq \epsilon+
  \sum_{n=(1-\delta) m(Y)N}^{(1+\delta) m(Y)N}
  m\{M(T_Y^n y) \geq a B_N\}.
  \end{equation*}

As $m$ is invariant by $T_Y$, we have $m\{M(T_Y^n y) \geq a B_N\}
=m\{ M \geq a B_N\} \leq \frac{C}{N}$. Thus,
  \begin{equation*}
  m\left\{ y \tq M\left(T_Y^{n(y,N)}y\right) \geq a B_N \right\}
  \leq \epsilon+
  2\delta m(Y)N \frac{C}{N}
  =2\epsilon.
  \end{equation*}
\qedfact
\end{proof}

\emph{Fact: $\frac{1}{B_N} \left(S^Y_{n(y,N)} f_Y- S^Y_{N
m(Y)}f_Y\right)$ tends to $0$ in probability on $X$ when $N \to
\infty$.}
\begin{proof}
By Lemma \ref{proba_sur_Y_donne_proba_sur_X}, it is sufficient to
prove it on $Y$. Without loss of generality, we can use the
natural extension and assume that $T_Y$ is invertible.

For $n<0$, write $S^Y_n f_Y=\sum_0^{|n|-1} f_Y\circ T_Y^{-j}$.
Then, setting $\nu(y,N)=n(y,N)-N m(Y)$,
\begin{equation}
  \label{definit_nu}
  S^Y_{n(y,N)} f_Y(y)- S^Y_{N m(Y)}f_Y(y)
  =S^Y_{\nu(y,N)} f_Y \left(T^{N m(Y)}(y)\right).
  \end{equation}
If $A>0$ and $N\in \N$, as $E_Y(\phi)=1/m(Y)$, we get
   \begin{multline*}
   \{y\tq  \nu(y,N) \geq A B_N^{1/b}\}
   =\{ n(y,N) \geq A B_N^{1/b}+ N m(Y)\}
   =\{ S^Y_{A B_N^{1/b}+Nm(Y)}\phi < N\}
   \\=
   \left \{ \frac{S^Y_{A B_N^{1/b}+N m(Y)}\phi - (A B_N^{1/b}+
   Nm(Y))E_Y(\phi)}{\left(B_{AB_N^{1/b}+Nm(Y)}\right)^{1/b}} <
   -\frac{A}{m(Y)}\left(\frac{B_N}{B_{AB_N^{1/b}+Nm(Y)}}\right)^{1/b}
   \right\}.
   \end{multline*}
For some integer $k$, we have $N \leq 2^k Nm(Y) \leq 2^k
(AB_{Nm(Y)}^{1/b}+Nm(Y))$. The assumption $\sup_{r\leq 2n}
\frac{B_r}{B_n} \leq C<\infty$ thus yields that
$\frac{B_N}{B_{AB_{Nm(Y)}^{1/b}+Nm(Y)}} \leq C^k$. In particular,
  \begin{equation*}
  \{y\tq  \nu(y,N) \geq A B_N^{1/b}\}
  \subset \left \{ \frac{S^Y_{A B_N^{1/b}+N m(Y)}\phi - (A B_N^{1/b}+
   Nm(Y))E_Y(\phi)}{\left(B_{AB_N^{1/b}+Nm(Y)}\right)^{1/b}} <
   -\frac{AC^{k/b}}{m(Y)}
   \right\}.
  \end{equation*}
Consequently, if $A$ is large enough, Assumption \ref{hypothese_4}
yields that $m\{y\tq  \nu(y,N) \geq A B_N^{1/b}\} \leq \epsilon$
for large enough $N$. We handle in the same way the set of points
where $\nu(y,N) \leq -A B_N^{1/b}$, using the assumption
$\inf_{r\geq n}\frac{B_r}{B_n}>0$.
%
%
%
%
We have thus proved:
  \begin{equation}
  \label{convergence_nuy}
  \forall \epsilon>0, \exists A>0, \exists N_0>0, \forall N \geq
  N_0,\
  m\{ y\tq |\nu(y,N)| \geq A B_N^{1/b} \} \leq \epsilon.
  \end{equation}

Set $W_N(y)=\frac{1}{B_N} S_{\nu(y,N)} f_Y \left(T_Y^{ N
m(Y)}(y)\right)$, we will show that it tends to $0$ in
distribution, which will conclude the proof, by
\eqref{definit_nu}. Take $a>0$, we show that $m(|W_N|>a) \to 0$
when $N\to \infty$.

Let $\epsilon>0$. Assumption \ref{hypothese_3} ensures that there
exists $\tilde{Y}$ with $m(\tilde{Y})\geq m(Y)-\epsilon$ and $N_1$
such that $\frac{1}{|N|^b}|S_N^Y f_Y| \leq \epsilon$ on
$\tilde{Y}$, for every $|N|\geq N_1$. Define $Y'_N=\{y\in Y \tq
|\nu(y,N)|< N_1\}$ and $Y''_N=\{y \in Y \tq |\nu(y,N)|\geq N_1\}$.
We estimate first the contribution of $Y'_N$.

Set $\psi(y)=\sum_{-N_1}^{N_1-1} |f_Y\circ T_Y^j|$. Since $\psi$
is measurable, there exists a constant $C$ and a subset $Z$ of $Y$
with $m(Z)\geq m(Y)-\epsilon$ and $\psi \leq C$ on $Z$. Then, for
$y\in Y'_N$, we have $|W_N(y)|\leq \frac{1}{B_N} \psi\left(T_Y^{N
m(Y)} y\right)$. Set $Z_N=Y'_N \cap T_Y^{- N m(Y)}(Z)$: it
satisfies $m(Z_N) \geq m(Y'_N)-\epsilon$. On $Z_N$, we have
$|W_N|\leq \frac{C}{B_N}$, whence, for large enough $N$, $|W_N|<a$
on $Z_N$. Thus, for large enough $N$,
  \begin{equation*}
  m\left\{y \in Y'_N \tq |W_N(y)|\geq a\right\}
  \leq m\left\{y\in Z_N \tq |W_N(y)|\geq a\right\}+\epsilon
  =\epsilon.
  \end{equation*}

We estimate then the contribution of $Y''_N$. Set
$\tilde{Y}''_N=Y''_N \cap T_Y^{-N m(Y)}(\tilde{Y})$, satisfying
$m(\tilde{Y}''_N) \geq m(Y''_N)-\epsilon$. Thus,
  \begin{equation*}
  m(|W_N|\geq a) \leq m\{y\in \tilde{Y}''_N \tq |W_N(y)|\geq
  a\}+2\epsilon.
  \end{equation*}
On $\tilde{Y}''_N$, $|\nu(y,N)|\geq N_1$, whence
$\frac{1}{|\nu(y,N)|^b}\left|S^Y_{\nu(y,N)}f_Y \left(T_Y^{N
m(Y)}y\right)\right| \leq \epsilon$. Thus, $|W_N(y)|\leq \epsilon
\frac{|\nu(y,N)|^b}{B_N}=\epsilon
\left(\frac{|\nu(y,N)|}{B_N^{1/b}}\right)^b$. Consequently,
  \begin{equation*}
  m(|W_N|\geq a) \leq m\left (\frac{|\nu(y,N)|}{B_N^{1/b}} \geq
  \left(\frac{a}{\epsilon}\right)^{1/b} \right) +2\epsilon.
  \end{equation*}
Thus, if $\epsilon$ is small enough, and $N$ large enough,
\eqref{convergence_nuy} yields that $m(|W_N|\geq a) \leq
3\epsilon$. \qedfact
\end{proof}

The three facts we have just proved imply that $\frac{S_N
f(y,i)}{B_N}-\frac{S_{Nm(Y)}^Y f_Y(y)}{B_N} \to 0$ in distribution
on $X$. As $\frac{S_{Nm(Y)}^Y f_Y(y)}{B_N} \to Z$ in distribution
on $X$, by \eqref{Thm_limite_SY}, this concludes the proof.
\end{proof}

\section{Multiple decorrelations and $L^p$-boundedness}

\label{appendice:pene}

The following theorem has been useful in this paper:

\begin{thm}
\label{thm_borne_Lp_pene} Let $F: \omega \to 4\omega$ on the
circle $S^1$. Then, for every $p \in [1,\infty)$, there exists a
constant $K_p$ such that, for every $n\in \N$, for every
$f_0,\ldots,f_{n-1}:S^1 \to \R$ bounded by $1$, of zero average and
$1$-Lipschitz,
  \begin{equation*}
  \norm{\sum_{k=0}^{n-1} f_k \circ F^k }_p \leq K_p \sqrt{n}.
  \end{equation*}
\end{thm}

This result has essentially been proved by Françoise Pène, in a
much broader context. Her proof depends on a property of multiple
decorrelations, which is implied by the spectral gap of the
transfer operator:

\begin{lem}
Let $\norm{f}$ be the Lipschitz norm of the function $f$ on the
circle $S^1$. Then, for every $m,m'\in \N$, there exist $C>0$ and
$\delta<1$ such that, for every $N\in \N$, for every increasing
sequences $(k_1,\ldots,k_m)$ and $(l_1,\ldots,l_{m'})$, for every
Lipschitz functions $G_1,\ldots,G_m,H_1,\ldots, H_{m'}$,
  \begin{equation}
  \label{decorr_mult}
  \left|\Cov\left(\prod_{i=1}^m G_i \circ F^{k_i}, \prod_{j=1}^{m'}
    H_j\circ F^{N+l_j}\right)\right|
  \leq C \left(\prod_{i=1}^m \norm{G_i}\right) \left(\prod_{j=1}^{m'}
    \norm{H_j} \right) \delta^{N-k_m}.
  \end{equation}
\end{lem}
Here $\Cov(u,v)=\int uv -\int u \int v$.

\begin{proof}
Let $\hat{F}$ be the transfer operator associated to $F$, and
acting on Lipschitz functions. It is known that it admits a
spectral gap and that its iterates are bounded, i.e.\ there exist
constants $M>0$ and $\delta<1$ such that $\bigl\|\hat{F}^n
f\bigr\|\leq M \norm{f}$, and $\bigl\|\hat{F}^n f\bigr\|\leq
M\delta^n \norm{f}$ if $\int f=0$.

We can assume that $N\geq k_m$ (otherwise, $\delta^{N-k_m}\geq 1$,
and the inequality \eqref{decorr_mult} becomes trivial). Then,
writing $\phi=\prod_{i=1}^m G_i \circ F^{k_i}$ and
$\psi=\prod_{j=1}^{m'} H_j\circ F^{l_j}$, we get
  \begin{align*}
  \left|\Cov(\phi,\psi\circ F^N)\right|
  &
  =\left|
  \int \left(\phi-\int \phi\right) \psi\circ F^N\right|
  =\left|\int \hat{F}^N \left(\phi-\int \phi\right) \psi \right|
  \\&
  \leq \norm{\hat{F}^N \left(\phi-\int \phi\right)} \norm{\psi}_\infty.
  \end{align*}
But
  \begin{align*}
  \hat{F}^N(\phi)&
  =\hat{F}^N \left(\prod G_i^{k_i}\right)
  =\hat{F}^{N-k_m} (G_m \hat{F}^{k_m-k_{m-1}}(G_{m-1}
  \hat{F}^{k_{m-1}-k_{m-2}}(\ldots \hat{F}^{k_2-k_1} (G_1))\ldots)
  \\&
  =:\hat{F}^{N-k_m}(\chi).
  \end{align*}
As the iterates of $\hat{F}$ are bounded on Lipschitz functions,
we get a bound on the Lipschitz norm of $\chi$: $\norm{\chi}\leq
M^{m-1} \prod \norm{G_i}$. Moreover, $\int \chi=\int \phi$, whence
  \begin{align*}
  \norm{\hat{F}^N \left(\phi-\int \phi\right)}
  &
  =\norm{\hat{F}^{N-k_m} \left(\chi -\int \chi\right)}
  \leq M\delta^{N-k_m} \norm{\chi-\int \chi}
  \\&
  \leq M\delta^{N-k_m} M^{m-1} \prod \norm{G_i}.
  \end{align*}
\end{proof}

When $p$ is an even integer, Theorem \ref{thm_borne_Lp_pene} is
then a consequence of \cite[Lemma 2.3.4]{pene:averaging}. The
Hölder inequality gives the general case.

\begin{rmq}
\label{remarque_pene}
The same result holds for Hölder functions instead of Lipschitz
functions, with the same proof.
\end{rmq}

We will also need the following result:
\begin{thm}
\label{ineg_max_abstraite}
Let $T$ be a measure preserving transformation on a space $X$. Let
$f:X\to \R$ and $p>2$ be such that
  \begin{equation*}
  \exists C>0, \forall n\in \N^*, \norm{S_n f}_p \leq C \sqrt{n}.
  \end{equation*}
Write $M_nf(x)=\sup_{1\leq k\leq n} |S_k f(x)|$. Then there exists
a constant $K$ such that
  \begin{equation*}
  \forall n\geq 2, \norm{M_n f}_p \leq K (\ln n)^{\frac{p-1}{p}}
  \sqrt{n}.
  \end{equation*}
\end{thm}
\begin{proof}
Let $n\in \N^*$. Let $k<2^n$, and write its binary decomposition
$k=\sum_{j=0}^{n-1} \epsilon_j 2^j$, with $\epsilon_j \in
\{0,1\}$. Set $q_j=\sum_{l=j}^{n-1} \epsilon_l 2^l$ (in
particular, $q_0=k$ and $q_n=0$). Then $S_k f=\sum_{j=0}^{n-1}
(S_{q_j}f-S_{q_{j+1}}f)$. Consequently, the convexity inequality
$(a_0+\ldots +a_{n-1})^p \leq n^{p-1} (a_0^p+\ldots+a_{n-1})^p$
gives that
  \begin{equation*}
  |S_k f|^p \leq n^{p-1} \sum_{j=0}^{n-1}
  |S_{q_j}f-S_{q_{j+1}}f|^p.
  \end{equation*}
Note that $q_{j+1}$ is of the form $\lambda 2^{j+1}$ with $0\leq
\lambda \leq 2^{n-j-1}-1$, and $q_j$ is equal to $q_{j+1}$ or
$q_{j+1}+2^j$. Thus,
  \begin{equation*}
  |S_k f|^p \leq n^{p-1} \sum_{j=0}^{n-1} \left(
  \sum_{\lambda=0}^{2^{n-j-1} -1} \left|S_{\lambda 2^{j+1} +2^j}f
  -S_{\lambda 2^{j+1}}f \right|^p \right).
  \end{equation*}
The right hand term is independent of $k$, and gives a bound on
$|M_{2^n -1} f|^p$. Moreover,
  \begin{equation*}
  \int \left|S_{\lambda 2^{j+1} +2^j}f -S_{\lambda 2^{j+1}}f
  \right|^p =\int \left|S_{2^j} f\right|^p \leq C^p \sqrt{2^j}^{p}.
  \end{equation*}
Thus, we get
  \begin{equation*}
  \int |M_{2^n -1} f|^p \leq
  n^{p-1} \sum_{j=0}^{n-1} 2^{n-j-1} C^p 2^{pj/2}
  \leq K n^{p-1} 2^n 2^{(\frac{p}{2}-1)n}
  = K n^{p-1} \sqrt{2^n}^p.
  \end{equation*}
For times of the form $2^n-1$, this is a bound of the form
$\norm{M_t}_p \leq K (\ln t)^{\frac{p-1}{p}} \sqrt{t}$. To get the
same estimate for an arbitrary time $t$, it is sufficient to
choose $n$ with $2^{n-1} \leq t<2^n$, and to note that $M_t \leq
M_{2^n -1}$.
\end{proof}

\begin{cor}
\label{ineg_max}
Let $F: \omega \to 4\omega$ on the circle $S^1$, let $\chi:S^1 \to
\R$ be a Hölder function with $0$ average, and $p>2$. Write $M_n
\chi(x)=\sup_{1\leq k\leq n} |S_k \chi(x)|$. Then there exists a
constant $K$ such that
  \begin{equation*}
  \forall n\geq 2, \norm{M_n \chi}_p \leq K (\ln
  n)^{\frac{p-1}{p}} \sqrt{n}.
  \end{equation*}
\end{cor}
\begin{proof}
Theorem \ref{thm_borne_Lp_pene} (or rather the remark following
it, for the Hölder case) shows that $\norm{S_n \chi}_p \leq
C\sqrt{n}$. Consequently, Theorem \ref{ineg_max_abstraite} gives
the conclusion.
\end{proof}

\bibliography{biblio}

\begin{thebibliography}{dMvS93}

\bibitem[Aar97]{aaronson:book}
Jon Aaronson.
\newblock {\em An introduction to infinite ergodic theory}, volume~50 of {\em
  Mathematical Surveys and Monographs}.
\newblock American Mathematical Society, 1997.

\bibitem[AD98]{aaronson_denker:central}
Jon Aaronson and Manfred Denker.
\newblock A local limit theorem for stationary processes in the domain of
  attraction of a normal distribution.
\newblock Preprint, 1998.

\bibitem[AD01]{aaronson_denker}
Jon Aaronson and Manfred Denker.
\newblock Local limit theorems for partial sums of stationary sequences
  generated by {G}ibbs-{M}arkov maps.
\newblock {\em Stoch. Dyn.}, 1:193--237, 2001.

\bibitem[dMvS93]{demelo_vanstrien}
Welington de~Melo and Sebastian van Strien.
\newblock {\em One-dimensional dynamics}, volume~25 of {\em Ergebnisse der
  Mathematik und ihrer Grenzgebiete : 3}.
\newblock Springer, 1993.

\bibitem[Fel66]{feller:2}
William Feller.
\newblock {\em An Introduction to Probability Theory and its Applications,
  volume 2}.
\newblock Wiley Series in Probability and Mathematical Statistics. J. Wiley,
  1966.

\bibitem[GH88]{guivarch-hardy}
Yves Guivarc'h and Jean Hardy.
\newblock Théorèmes limites pour une classe de chaînes de {M}arkov et
  applications aux difféomorphismes d'{A}nosov.
\newblock {\em Ann. Inst. H. Poincaré Probab. Statist.}, 24:73--98, 1988.

\bibitem[Gou02]{gouezel:stable}
S{\'e}bastien Gou{\"e}zel.
\newblock Central limit theorem and stable laws for intermittent maps.
\newblock Preprint, 2002.

\bibitem[Hu01]{hu:almost_hyperbolic}
Huyi Hu.
\newblock Statistical properties of some almost hyperbolic systems.
\newblock In {\em Smooth ergodic theory and its applications (Seattle, WA,
  1999)}, volume~69, pages 367--384. Amer. Math. Soc., 2001.

\bibitem[Kac96]{vitesse_birkhoff}
A.~G. Kachurovski\u\i.
\newblock Rates of convergence in ergodic theorems.
\newblock {\em Russian Math. Surveys}, 51:653--703, 1996.

\bibitem[LSV99]{liverani_saussol_vaienti}
Carlangelo Liverani, Beno{\^\i}t Saussol, and Sandro Vaienti.
\newblock A probabilistic approach to intermittency.
\newblock {\em Ergodic Theory and Dynamical Systems}, 19:671--685, 1999.

\bibitem[MT02]{melbourne_torok}
Ian Melbourne and Andrew T{\"o}r{\"o}k.
\newblock Statistical limit theorems for suspension flows.
\newblock Preprint, 2002.

\bibitem[P{\`e}n02]{pene:averaging}
Fran{\c{c}}oise P{\`e}ne.
\newblock Averaging method for differential equations perturbed by dynamical
  systems.
\newblock {\em ESAIM Probab. Statist.}, 6:33--88 (electronic), 2002.

\bibitem[PY01]{pollicott_yuri:indifferent}
Mark Pollicott and Michiko Yuri.
\newblock Statistical properties of maps with indifferent periodic points.
\newblock {\em Comm. Math. Phys.}, 217(3):503--520, 2001.

\bibitem[Via97]{viana:multidim_attr}
Marcelo Viana.
\newblock Multidimensional nonhyperbolic attractors.
\newblock {\em Publ. Math. IHES}, 85:63--96, 1997.

\end{thebibliography}
\bibliographystyle{alpha}

\end{document}